\documentclass{article}
\usepackage{amsfonts, amssymb, latexsym, graphicx}
\usepackage[all]{xy}
\input xy
\xyoption{all}
\newtheorem{definition}{Definition}
\newtheorem{theorem}{Theorem}
\newtheorem{proposition}{Proposition}

\begin{document}
\title{Pencils of cubics with eight base points lying in convex position in $\mathbb{R}P^2$}
\author{S\'everine Fiedler-Le Touz\'e}
\maketitle
\begin{abstract}
To a generic configuration of eight points in convex position in the real plane, we associate a {\em list\/} consisting of the following information: for all of the $56$ conics determined by five of the points, we specify the position of each of the three remaining points, inside or outside. We prove that the number of possible lists, up to the action of $D_8$ on the set of points, is $47$, and we give two possible ways of encoding these lists. 
A generic complex pencil of cubics has twelve singular (nodal) cubics and nine distinct base points, any eight of them determine the ninth one, hence the pencil. If the base points are real, exactly eight of these singular cubics are distinguished, that is to say real with a loop containing some base points. We call combinatorial cubic a topological type (cubic, base points), and combinatorial pencil the sequence of eight successive combinatorial distinguished cubics.
Let us choose representants of the $47$ orbits.  In most cases, but four exceptions, a list determines a unique combinatorial pencil. We get a complete classification of the combinatorial pencils of cubics with eight base points in convex position. 
Up to the action of $D_8$ on the set of points, there are $43$ such pencils.
\end{abstract}

\section{Pencils of cubics} \label{sec:1}

\subsection{Preliminaries} \label{subsec: 1.1}

Given nine generic points in $\mathbb{C}P^2$, there exists one single cubic
passing through them. Given eight generic points in $\mathbb{C}P^2$, there 
exists a one-parameter family of cubics passing through them. We will call 
such a family a {\em pencil of cubics\/}. Let $F_0$ and $F_1$ be two cubics 
of a pencil $\mathcal{P}$. They intersect at a ninth point. As the other 
cubics of $\mathcal{P}$ are linear combinations of $F_0$ and $F_1$, they all 
pass through this ninth point. We call these nine points the 
{\em base points\/} of $\mathcal{P}$. If eight of the base points are real, the
pencil is real, and hence the ninth base point is also real.
A pencil of cubics is a line in the space $\mathbb{C}P^9$ of complex cubics. 
Let $\Delta$ be the discriminantal hypersurface of $\mathbb{C}P^9$, formed by 
the singular cubics. The hypersurface $\Delta$ is of degree 12. 
Hence, a generic pencil of cubics intersects $\Delta$ transversally at 12
regular points.
Otherwise stated, a generic pencil $\mathcal{P}$ has exactly
12 singular (nodal) cubics. A non-generic pencil will be called {\em singular
pencil\/}.
Let $\mathcal{P}$ be a real pencil with nine real base points and denote 
the real part of $\mathcal{P}$ by $\mathbb{R}\mathcal{P}$. 
Let $n \leq 12$ be the number of real singular cubics of $\mathcal{P}$. 
Let $C_3$ be one of these cubics. The double point $P$ of $C_3$ is 
{\em isolated\/} if the tangents to $C_3$ at $P$ are non-real, otherwise $P$ 
is {\em non-isolated\/}. If $P$ is non-isolated, $C_3 \setminus P= 
\mathcal{J} \cup \mathcal{O}$, where
$[\mathcal{J} \cup P] \not= 0$ and $[\mathcal{O} \cup P]=0$ in 
$H_1(\mathbb{R}P^2)$.
We say that $\mathcal{O}$ is the loop and $\mathcal{J}$ is the odd component 
of $C_3$. Notice that the loop $\mathcal{O}$ is convex.
The estimation of $n$ presented hereafter is due to V.Kharlamov, see
\cite{D-K}.

One has: $n = n_1+n_2+n_3$, where $n_1$ is the number of cubics with an
isolated double point, $n_2$ is the number of cubics with a loop containing 
no base points, and $n_3$ is the number of cubics with a loop containing some 
base points.To evaluate $n$, one recalculates the Euler characteristic of 
$\mathbb{R}P^2$, fibering $\mathbb{R}P^2$ with the cubics of $\mathbb{R}\mathcal{P}$. 
Each isolated double point, and each base point contributes by $+1$; and 
each non-isolated double point contributes by $-1$, so that one gets:

\begin{displaymath}
1=\chi(\mathbb{R}P^2)=9+n_1-(n-n_1)
\end{displaymath}

So, $n-2n_1=8$. Thus, $n=8, 10$ or $12$ and correspondingly, $n_1=0, 1$
or $2$.

Consider a motion in $\mathbb{R}\mathcal{P}$ starting from a cubic with an
isolated double point. If we choose the direction of the motion properly, an 
oval appears, grows, and attaches itself to the odd component, forming a loop 
that contains no base point.
Conversely, starting from a cubic with a loop containing no base point, one 
can move in the pencil so that there appears a cubic with an oval. As this 
oval lies inside of the loop, it shrinks when one moves further, and 
degenerates into an isolated double point.
Thus, $n_2=n_1$ and $n_3=8$ independently of $n$. Let us call the eight
cubics of the third type the {\em distinguished cubics\/} of $\mathcal{P}$.
We shall picture $\mathbb{R}\mathcal{P}$ by a circle, divided in eight 
portions by the eight distinguished cubics. Let us remark that this number 
$n_3=8$ is the Welschinger invariant $W_3$, see \cite{W}. The number of real rational plane curves of degree $d$ going through $3d-1$ generic 
points of $\mathbb{R}P^2$ is always finite. Let $c_1$ be the number of such 
curves with an even number of isolated nodes, and $c_2$ be the number of 
such curves with an odd number of isolated nodes. Welschinger proved that 
the difference $W_d=c_1-c_2$ does not depend on the
choice of the $3d-1$ points. 
Let $1, \dots 8$ be eight generic points in $\mathbb{R}P^2$, determining a pencil of cubics $\mathcal{P}$. Let $9$ be the ninth base point, and $C_3$ be a real cubic of $\mathcal{P}$. We call {\em combinatorial cubic\/} $C_3$ the topological type of $(C_3, 1, \dots, 9)$, and {\em combinatorial pencil\/} the cyclic sequence of the successive eight combinatorial distinguished cubics.
From now on, the same notation $\mathcal{P} =\mathcal{P}(1, \dots 8)$ will refer to the (real part of the) pencil of cubics determined by $1, \dots 8$, and to the combinatorial pencil. It will always be clear from the context whether we speak of a combinatorial or of an actual pencil.

\begin{theorem}
Up to the action of $D_8$, eight generic points $1, \dots 8$ lying in convex position in $\mathbb{R}P^2$ may realize $43$ combinatorial pencils $\mathcal{P}(1, \dots 8)$. 
\end{theorem}

\subsection{Singular pencils} \label{subsec: 1.2}

Let $\mathcal{P}$ be an (actual) pencil of cubics with only real base points.
Move these points till $\mathcal{P}$
degenerates into a singular pencil $\mathcal{P}_{sing}$. The degeneration is 
generic if and only if $\mathcal{P}_{sing}$ intersects $\Delta$ transversally
at 10 regular points and
\begin{enumerate}
\item
$\mathcal{P}_{sing}$ is tangent to $\Delta$ at one regular point, or
\item
$\mathcal{P}_{sing}$ crosses transversally a stratum of codimension 1 of 
$\Delta$
\end{enumerate}

\begin{figure}
\centering
\includegraphics{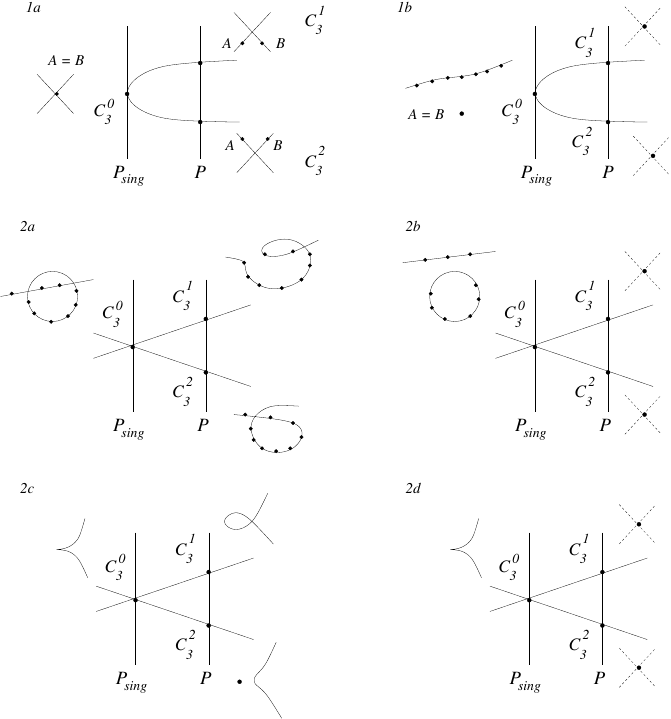}
\caption{\label{delta} The degeneration of $\mathcal{P}$ into $\mathcal{P}_{sing}$} 
\end{figure}

In both cases, two singular cubics $C^1_3$ and $C^2_3$ of $\mathcal{P}$ come 
together to yield one singular cubic $C_3^0$ of $\mathcal{P}_{sing}$. The cubic 
$C_3^0$ is necessarily real, the cubics  $C^1_3$ and $C^2_3$ are either both real 
or complex conjugated. 
Let $F$ be a generic cubic of $\Delta$ with node at some point $p$,
then $T_F\Delta=\{F+G \vert  G(p)=0\}$. Therefore, $\mathcal{P}_{sing}$ 
satisfies the condition 1) if and only if $\mathcal{P}_{sing}$ has a double 
base point, at $p$. Otherwise stated, $\mathcal{P}_{sing}$ is
obtained from $\mathcal{P}$ by letting two base points $A$ and $B$ of 
$\mathcal{P}$ come together. Move $A$ towards $B$ along the line $(AB)$.
For simplicity, we assume that $B$ is the origin $(0, 0)$ of the plane, and the 
direction is the $x$-axis. The condition that some cubic  $H$ 
passes through $A$ and $B$ becomes at the limit: $H(0, 0)=0$ and 
$\partial{H} \over{\partial{x}}$$(0, 0) =0$. All of the cubics of $\mathcal{P}_{sing}$ but one are tangent to the $x$-axis. A unique cubic $C_3^0$ is singular at $A = B$ (for this cubic, the other partial derivative $\partial{H} \over{\partial{y}}$$(0, 0)$ is also equal to $0$).
In case 2), the cubic $C_3^0$ must be reducible (product of a line and a conic), or 
have a cusp. If $C_3^0$ is reducible, the genericity imposes that three of the
base points lie on the line, and the other six lie on the conic.
If $C_3^0$ has a cusp, the cubics $C^1_3$ and  $C^2_3$ are either non-distinguished real cubics, one with a loop, the  other with an isolated double point, or 
complex conjugated. See Figure~\ref{delta}, where the various possible types for the degeneration $\mathcal{P} \to \mathcal{P}_{sing}$ have been denoted by 1a, 1b, 2a, 2b, 2c, 2d. 
In the right-hand part of the figure, the dotted crosses symbolize the non-real
nodes of the cubics $C^1_3$ and  $C^2_3$. Moving further, one gets a new pencil
$\mathcal{P}'$. The cubic $C_3^0$ is replaced by a new pair of nodal cubics $C_3^3$, $C_3^4$. If $\mathcal{P} \to \mathcal{P}_{sing}$ realizes one of the types 1a, 1b, 2a or 2b, $\mathcal{P}' \to \mathcal{P}_{sing}$ is of the same type, whereas the types 2c and 2d are swapped.  

Let us describe more precisely the case 1).
We call {\em elementary arc\/} $AB$ an arc connecting $A$ to $B$ and containing 
no other base point. As $\mathcal{P}$ can degenerate into $\mathcal{P}_{sing}$
letting $A$ and $B$ come together, some cubics of $\mathcal{P}$ must have an 
elementary arc $AB$. Start from such a cubic and move in any direction in 
$\mathbb{R}\mathcal{P}$. At some moment, the mobile arc $AB$ must glue to 
another arc, and then disappear. Thus $\mathcal{P}$ has two distinguished
cubics corresponding to the openings of the arc $AB$. We call {\em singular 
elementary arc $AB$\/} the non-smooth arc $AB$ of either of these cubics.
For 1a, these two cubics are  $C^1_3$ and $C^2_3$, they come together to 
yield the cubic $C_3^0$, which has a non-isolated double point at $A=B$. 
All  three combinatorial cubics $C_3^0$,  $C^1_3$ and $C^2_3$ are identical 
outside of a neighbourhood of $A \cup B$. Note that the pair of cubics 
$C^1_3$, $C^2_3$ divide the pencil $\mathbb{R}\mathcal{P}$ in two portions, in one portion, the cubics have an elementary arc $AB$, in the other they have none. The latter portion disppears in the motion $\mathcal{P} \to
\mathcal{P}_{sing}$. For 1b, the two cubics corresponding to the opening of
$AB$ divide the pencil in two portions, the cubics in each portion have respectively one and two elementay arc(s) $AB$. The cubics of the latter portion have all an oval passing through $A, B$, and no other base point. When the motion reaches
$\mathcal{P}_{sing}$, a nodal cubic $C_3^0$ appears inside of this portion.
Percoursing this portion, one sees the oval shrink till it becomes an isolated node. Then the oval reappers on the other side and grows.  


Let us finally mention that eight points in $\mathbb{R}P^2$ determine a pencil of cubics, with a ninth base point, if and only if no four of the eight points are aligned, and no seven are coconic. If three of the points are aligned, the ninth base point lies on the conic determined by the five others; if six points are coconic, the ninth base point is on the line determined by the other two: the pencil is singular. If four points are aligned, one gets a pencil of reducible cubics, product of the line and the pencil of conics determined by the other four points. If seven points are coconic, one gets again a pencil of reducible cubics, product of the conic and the pencil of lines through the remaining point. 
If two points coincide, say $A$ and $B$, one gets a one-parameter family of
singular pencils: consider the cubic $C_3^0$ through $1, \dots 8$, with node at $A = B$. The ninth point $9$ may be chosen anywhere on this cubic. Finally, six generic points $1, \dots 6$ in $\mathbb{R}P^2$ determine six rational pencils, with respective nodes at $1, \dots 6$. Each pencil has five reducible cubics. Let us call {\em combinatorial cubic $C_3$\/} a topological type $(C_3, 1, \dots 6)$ and {\em combinatorial pencil\/} the cyclic sequence of five combinatorial reducible cubics. Up to the action of $\mathcal{S}_6$ on $1, \dots 6$, there are four lists of six combinatorial pencils, \cite{F3}.  

\section{Configurations of eight points lying in convex position in $\mathbb{R}P^2$, lists of conics}

\subsection{Mutual position of points and conics}

We say that $n$ points lie in convex position in
$\mathbb{R}P^2$ if there exists a line $L \subset \mathbb{R}P^2$ such that
the $n$ points lie in convex position in the affine plane 
$\mathbb{R}P^2 \setminus L$. Note that any set of $n \leq 5$ points in
$\mathbb{R}P^2$ lie in convex position.
Consider $n \geq 5$ generic points lying in strictly convex position in
$\mathbb{R}P^2$, say $1, \dots, n$. Let $L(1, \dots, n)$ be the list of the 
${ n \choose 5 }$ conics through five of these points, enhanced for each conic, with 
the position of each of the remaining $n-5$ points (inside or outside).
So a (generic) list, written in full extent, consists of $(n-5)\times { n \choose 5 }$ {\em elementary data\/} of the type $k < C_2$ ($k$ lies inside of the conic $C_2$) or $k > C_2$ ($k$ lies outside of $C_2$). In section 4, we will consider also non-generic configuration of points $1, \dots, 8$ with six of them on a conic, and the corresponding non-generic lists

\begin{quote}
{\em How many different possibilities can be realized by $L(1, \dots, n)$ 
when one lets the points $1, \dots, n$ move?\/}
\end{quote}

For $n=6$, note that the condition $6 > 12345$ is equivalent to $5 < 12346$.
One single elementary data determines $L(1, \dots 6)$. There are two lists:
$6 < 12345$ and $6 > 12345$. We write shortly: $\natural L(1, \dots, 6)=2$.
The two lists, written in full extent are displayed in Figure~\ref{listsix}.
\begin{figure}
\begin{tabular}{|c|c|c|c|c|}
\hline
$L(1, \dots 6)$ & \multicolumn{2}{c|}{} &  \multicolumn{2}{c|}{}\\
\hline
$C_2$ & in & out & in & out\\
$12345$ & $6$ & & & $6$\\
$12346$ & & $5$ & $5$ & \\
$12356$ & $4$ & & & $4$\\
$12456$ & & $3$ & $3$ & \\
$13456$ & $2$ & & & $2$\\
$23456$ & & $1$ & $1$ & \\
\hline
\end{tabular}
\caption{\label{listsix}  The two lists $L(1, \dots 6)$} 
\end{figure}

Let now $n = 7$. Let $C_3$ be a cubic passing through the points $1, \dots, 7$. 
By abuse of language we will also call {\em cubic $C_3$\/} the topological type of 
$(C_3, 1, \dots, 7)$. Let $F$ be the equation of $C_3$ and $p$ be a point
among $1, \dots, 7$. The condition $F(p)=0$ is a linear equation in the
coefficients of $F$. If $p$ is a singular point of $C_3$, one gets
two supplementary linear equations:
$\partial{F} \over{\partial{x}}$$(p) =0$ and
$\partial{F} \over{\partial{y}}$$(p) =0$.Thus, as $1, \dots, 7$ are generic,
there exists exactly one real nodal cubic passing through $1, \dots, 7$ and 
having $p$ as double point.
Let $S(1, \dots, 7)$ be the list of the seven nodal cubics passing through the points $1, \dots, 7$, one of them being the double point. Denote by $C_3(k)$ the cubic
with node at $k$.  We may ask the same question as above, replacing $L(1, \dots 7)$ by $S(1, \dots 7)$.

\begin{proposition}
$\natural L(1, \dots, 7) = \natural S(1, \dots, 7) = 14$. 

The fourteen lists $S(1, \dots 7)$ are denoted by $1\pm, 2\pm, \dots, 7\pm$, they are all equivalent up to the action of the dihedral group $D_7$.

The combinatorial informations $S(1, \dots 7)$ and $L(1, \dots 7)$ are equivalent, and determined by only two elementary data: $S(1, \dots 7) = 1+$ if and only if $7 < 23456$ and $1 > 23456$.
\end{proposition}

\begin{figure}
\centering
\includegraphics{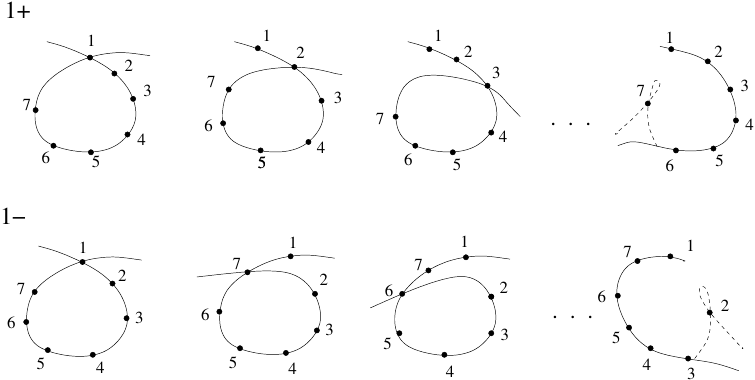}
\caption{\label{lists}  The lists $1+$ and $1-$}  
\end{figure}

The lists $1+$ and $1-$ are shown in Figure~\ref{lists}; the other lists 
$n\pm$ are obtained from $1\pm$ performing on $1, \dots, 7$ the cyclic 
permutation that replaces $1$ by $n$.
(The double point of the last cubic in each list may be an isolated node or 
a crossing. If it is a crossing, the loop is attached to the arc $61$ of the odd branch,
as shown in the Figure, or to the arc $12$.)
Note that any of the five non-extremal cubics of the list $S(1, \dots, 7)$ determines the whole of this list. Let us denote by $C_3(k)$ the cubic with node at $k$.

\begin{figure}
\centering
\includegraphics{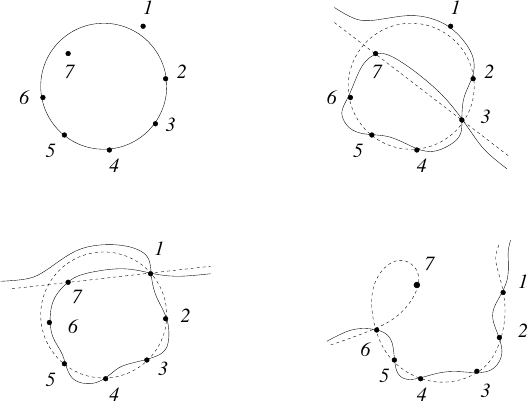}
\caption{\label{nodcub} Construction of the list $1+$} 
\end{figure}

\begin{description}
\item[Proof:]
The first six points $1, \dots 6$ may realize two different lists $L(1, \dots 6)$. For either of them, there are seven possible positions of the point $7$ with respect to the set of conics passing through five points among $1, \dots, 6$. One checks easily that the data {\em $L(1, \dots, 6)$, position of $7$\/} determines the list $L(1, \dots, 7)$. Thus, $\natural L(1, \dots, 7)=14$. 

Let $1, \dots 7$ be such that $7 < 23456$ and $1 > 23456$, see Figure~\ref{nodcub}. One may move $1$, preserving the convex position, towards $23456$ till $1$ reaches this conic.
The cubic $C_3(k)$ with $k = 2, 3, 4, 5$ or $6$ is obtained as perturbation of the reducible cubic $k7 \cup 23456$ letting $1$ move back to the outside of $23456$ (in Figure~\ref{nodcub}, we show the case $k = 3$). The condition {\em $1$ moves to the outside of $23456$\/} is equivalent to {\em
$6$ moves to the inside of $12345$\/}. The cubic $C_3(1)$ is obtained perturbing the reducible cubic $17 \cup 12345$. Assume that the cubic $C_3(7)$ has a loop, passing through some other of the seven points. This cubic cuts
then $C_3(6)$ at more than nine points, contradicting Bezout's theorem. Applying now Bezout's theorem with lines, we see that if $C_3(7)$ has a loop, it is attached to one of the arcs $61$ or $12$ of the pseudo-line. We have constructed 
$1+$, cyclic permutations and symmetries allow to obtain thirteen other lists.
As of now, we have proved that: $\natural L(1, \dots 7) = 14$, two elementary data of $L(1, \dots 7)$ determine $S(1, \dots 7)$, and $\natural S(1, \dots 7) \geq 14$.
Hence $\natural S(1, \dots 7) = 14$, the data $L(1, \dots 7)$ and $S(1, \dots 7)$
are equivalent.
The correspondences between the fourteen lists $S(1, \dots 7)$ and $L(1, \dots 7)$ (written in full extent) are displayed in Figures~\ref{lsevenone}-\ref{lseventwo}.
(For convenience, we have gathered some of the tabulars in the end.)
\end{description}

Another proof of this proposition, using cremona transformations, is presented in \cite{F2}. Let $A, B$ be two consecutive points in a configuration $1, \dots, 7$ realizing a generic list, and let $C_2$ be the conic through the other five points.
By definition, the {\em distance\/} $A \to B$ is $0$ if $A, B$ lie both inside or both outside of $C_2$; $+1$ if $A < C_2$ and $B > C_2$; $-1$ if $A > C_2$ and $B < C_2$. For the list $1+$, one has thus $1 \to 2 = \dots = 6 \to 7 = 0$, and $7 \to 1 = +1$.  

The first non trivial case is $n = 8$. It will turn out that the least number of elementary data necessary to determine a list $L(1, \dots 8)$ depends on the list. 
Let us denote by $\hat i$ the list determined by all of the points but $i$. 
The list $L(1, \dots, 8)$ may be encoded by the octuple of sublists $(\hat 1, \hat 2, \dots, \hat 8)$ with $\hat i =k\pm, k \in \{ 1, \dots, \hat i, \dots 8\}$. 
Consider the group $D_8$ of symmetries of the octagon, generated by  
the cyclic permutation $a = +1$ and the symmetry with respect to the 
axis $\sigma = 15$, see Figure~\ref{syms}. 

\begin{displaymath}
D_8=\{a, \sigma \vert a^8=id, \sigma^2=id , a\sigma=\sigma a^{-1}\}
\end{displaymath}

\begin{theorem}
Up to the action of $D_8$, eight generic points $1, \dots 8$ lying in convex position may realize $47$ different lists $L(1, \dots 8)$. The total number of generic lists is $752 = 47 \times 16$.
\end{theorem}

This theorem will be proved in two steps, in sections 2.2 (restriction part) and 2.6 (construction part).
\begin{figure}[htbp]
\begin{tabular}{ c c c c c c c c }
$15$ & 26 & 37 & 48 & (+1)(15) & (+1)(26) & (+1)(37) & (+1)(48)\\ 
$1 \leftrightarrow 1$ & $1 \leftrightarrow 3$ & $1 \leftrightarrow 5$ & $1 \leftrightarrow 7$ &  
$1 \leftrightarrow 2$ & $1 \leftrightarrow 4$ & $1 \leftrightarrow 6$ & $1 \leftrightarrow 8$\\
$2 \leftrightarrow 8$ & $2 \leftrightarrow 2$ & $2 \leftrightarrow 4$ & $2 \leftrightarrow 6$ &
$3 \leftrightarrow 8$ & $2 \leftrightarrow 3$ & $2 \leftrightarrow 5$ & $2 \leftrightarrow 7$\\
$3 \leftrightarrow 7$ & $4 \leftrightarrow 8$ & $3 \leftrightarrow 3$ & $3 \leftrightarrow 5$ &
$4 \leftrightarrow 7$ & $5 \leftrightarrow 8$ & $3 \leftrightarrow 4$ & $3 \leftrightarrow 6$\\
$4 \leftrightarrow 6$ & $5 \leftrightarrow 7$ & $6 \leftrightarrow 8$ & $4 \leftrightarrow 4$ &
$5 \leftrightarrow 6$ & $6 \leftrightarrow 7$ & $7 \leftrightarrow 8$ & $4 \leftrightarrow 5$\\
$5 \leftrightarrow 5$ & $6 \leftrightarrow 6$ & $7 \leftrightarrow 7$ & $8 \leftrightarrow 8$ & & & & 
\end{tabular}
\caption{\label{syms}  Action of $D_8$}
\end{figure}

\subsection{Admissible lists}

\begin{figure}
\centering
\includegraphics{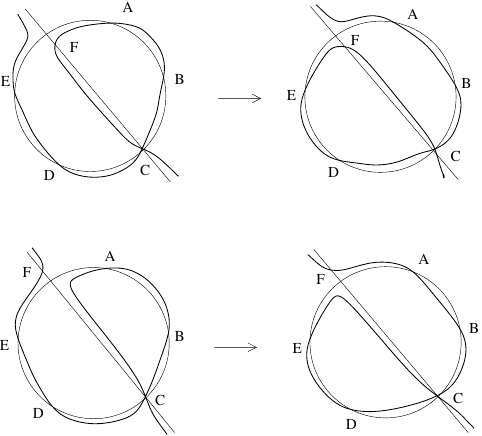}
\caption{\label{conics} Degenerations of the list $\hat G$} 
\end{figure}

Let $1, \dots 8$ lie in convex position and $\{ A, \dots G \} = \{ 1, \dots 7 \}$.
Move $8$, leaving the other points fixed and preserving the convex position. Consider an event of the form {\em $8$ crosses a conic $ABCDE$\/}, it induces a change of the list $\hat G$.  
The remaining point $F$ may lie inside or outside of the conic $ABCDE$, depending on the list $L(1, \dots, 7)$ (see upper and lower part of Figure~\ref{conics}). 
So there are: $21$ choices for the conic $ABCDE$, two choices of $G$,
and two possible positions of the last point $F$ with respect to $ABCDE$. 
Hence in total $84$ possibilities. 
Figure~\ref{conics} gathers all of these possibilities, showing how the cubic of
$\hat G$ with double point at $C$ changes when $8$ crosses $ABCDE$.
This cubic determines the whole list $\hat G$.
The point $8$ (not represented) may be placed in six different ways on
each cubic, according to the cyclic ordering of the set of points $\{8, A, B, 
C, D, E, F \}$. Assume for example that $8$ is situated between $A$ and $B$ 
in the cyclic ordering. When $8$ enters $ABCDE$, one has: 
$\hat G: E- \to A+$ (if $F$ inside of $ABCDE$), and $\hat G: F- \to
F+$ (if $F$ outside of $ABCDE$).
For each $G \in \{1, \dots, 7\}$, we find two possible chains of 
degenerations of $\hat G$ while moving the point $8$. The chain starting with $8\pm$ will be denoted by $\hat G\pm(8)$.
For each $\hat 8 \in \{ 1\pm, \dots 7\pm \}$ and each $G \in \{1, \dots, 7\}$ one watches which of the two possible chains is realized, see Figure~\ref{values}.
For each list $\hat 8 \in \{ 1 \pm, 2 \pm, 3 \pm, 4+ \}$, we draw a diagram whose rows are the chains of degenerations of $\hat 1, \dots \hat 7$. (We drop the other 
cases $\hat 8 \in \{7\pm, 6\pm, 5\pm, 4-\}$, that can be deduced from the 
first ones by the action of $D_8$).
Let $C_2$ and $C'_2$ be two adjacent conics in a column, we add a vertical arrow
from $C_2$ to $C'_2$ if the following holds:  {\em $8$ outside of $C_2$ implies $8$
outside of $C'_2$\/}.  See Figures~\ref{oneplus}-\ref{fourplus}.  

Chasing in the diagrams, we may find all of the admissible orbits, for the action of $D_8$, realizable by the lists $L(1, \dots 8)$. The explicit lists $L(1, \dots 8)$ obtained in this procedure are gathered in Figures~\ref{lone}-\ref{lsix}.
We denote these lists by $L_1, \dots L_{95}$, according to their appearance order.
Note that we cannot completely rule out redundancies: we get sometimes several representants of the same orbit. So, we choose for each orbit one representant that we write with normal fonts, the equivalent lists are written in bold. The lists written in normal fonts will be called for convenience {\em principal lists\/}, even if their choice is not canonical. 

Figures~\ref{lone}-\ref{ltwo} show the 64 admissible lists with $\hat 8 = 1+$. The first and the last are deduced one from the other by $\pm 1$. The 15 lists with $\hat 5 = 6+$ are mapped onto the 15 lists with $\hat 3 = 4+$ by $+3$.
The group of six lists with $\hat 2 = 1-$ splits into two subgroups that are mapped one onto the other by the symmetry $15$.
The group of 20 lists with $\hat 4 = 3-$ splits into two subgroups that are mapped one onto the other by $26$. 
The group of six lists with $\hat 6 = 5-$ splits into two subgroups that are mapped one onto the other by $37$.
Up to the action of $D_8$, there are 32 admissible lists with $\hat 8 = 1+$. 
Set now $\hat 8 = 1-$. First thing we rule out the orbit of $\hat 8 = 1+$, that is to say any list which is mapped by some element of $D_8$ onto a list with $\hat 8 = 1+$.
In other words, we set: $\hat 1 \not= 2+, 8-$, $\hat 2 \not= 3+, 1-$, $\hat 3 \not= 4+, 2-$, $\hat 4 \not= 5+, 3-$,  $\hat 5 \not= 6+, 4-$, $\hat 6 \not= 7+, 5-$ and  
$\hat 7 \not= 8+, 6-$.
Figure \ref{lthree} shows the six new admissible lists obtained. They split into two groups that are mapped one onto the other by $15$, so we are left with only three new orbits.  
For $\hat 8 = 2+$ we rule out the orbits of $\hat 8 = 1 \pm$, that is
we set: $\hat 1 \not= 2 \pm, 8 \pm$, $\hat 2 \not= 3 \pm, 1 \pm$, $\hat 3 \not= 4 \pm, 2 \pm$, $\hat 4 \not= 5 \pm, 3 \pm$,  
$\hat 5 \not= 6 \pm, 4 \pm$, $\hat 6 \not= 7 \pm, 5 \pm$ and  
$\hat 7 \not= 8 \pm, 6 \pm$. 
Figure \ref{lfour} shows the four new admissible lists with $\hat 8 = 2+$. They split into two groups that are mapped one onto the other by $15$, we get two new orbits.  
Set now $\hat 8 = 2-$. 
Once we have ruled out the orbits of $\hat 8 = 1 \pm, 2+$, we get the 13 new admissible lists shown in Figure~\ref{lfive}. Some of them are deduced from each other by $\pm 2, 15$ or $26$ so that there are only eight new orbits.
Let $\hat 8 = 3+$, once we have ruled out the orbits of $\hat 8 = 1 \pm, 2 \pm$,
we find  no new admissible lists.
For $\hat 8 = 3-$,  after excluding the orbits of $\hat 8 = 1 \pm, 2 \pm, 3+$,
we get the eight new admissible lists shown in Figure \ref{lsix}. They split into two groups that are mapped one onto the other by $\pm 2$, the first group splits in two subgroups that are mapped one onto the other by $15$. We get two new orbits.
At last, let $\hat 8 = 4+$, once we have excluded the orbits of $\hat 8 = 1 \pm, 2 \pm, 3 \pm$,
we find no new admissible lists.

The total number of admissible orbits is $47$. 
One checks easily that the 16 lists in each orbit are all distinct: assume that a principal list, say with $\hat 8 = 1+$, is invariant for some element of $D_8$, say $(+1)$. This list verifies: $\hat 8 = 1+$ and $\hat 1 = 2+$. There is no such list.
The total number of admissible lists is thus $47 \times 16 = 752$.
The non-generic lists invariant for some elements of $D_8$ will be studied in section 4.3.

\begin{figure}[htbp]
\begin{tabular}{ c | c }
 & $\hat 8$\\
\hline
$\hat 1 +$ & $1+, 2+, 2-, 4+, 4-, 6+, 6-$\\
$\hat 1 -$ & $1-, 3+, 3-, 5+, 5-, 7+, 7-$\\
$\hat 2 +$ & $1+, 1-, 2-, 4+, 4-, 6+, 6-$\\
$\hat 2 -$ & $2+, 3+, 3-, 5+, 5-, 7+, 7-$\\
$\hat 3 +$ & $1+, 1-, 3+, 4+, 4-, 6+, 6-$\\
$\hat 3 -$ & $2+, 2-, 3-, 5+, 5-, 7+, 7-$\\
$\hat 4 +$ & $1+, 1-, 3+, 3-, 4-, 6+, 6-$\\
$\hat 4 -$ & $2+, 2-, 4+, 5+, 5-, 7+, 7-$\\
$\hat 5 +$ & $1+, 1-, 3+, 3-, 5+, 6+, 6-$\\
$\hat 5 -$ & $2+, 2-, 4+, 4-, 5-, 7+, 7-$\\
$\hat 6 +$ & $1+, 1-, 3+, 3-, 5+, 5-, 6-$\\
$\hat 6 -$ & $2+, 2-, 4+, 4-, 6+, 7+, 7-$\\
$\hat 7 +$ & $1+, 1-, 3+, 3-, 5+, 5-, 7+$\\
$\hat 7 -$ & $2+, 2-, 4+, 4-, 6+, 6-, 7-$\\  
\end{tabular}
\caption{\label{values}
Possible values of $\hat 8$ for the chains $\hat n \pm = \hat n \pm(8), n =1, \dots, 7$}
\end{figure}

\subsection{Extremal lists}

Consider a configuration of eight points $1, \dots 8$ lying in convex position in the plane.
Any piece of information $\hat P = n\pm$, with $P \in \{1, \dots 8\}$  
is equivalent to a statement of the form $F < C_2$ {\em and\/} $G > C_2$, where $F, G$ are two points among $1, \dots, 8$, different from $P$.
For example, $\hat 8 = 1+$ if and only if $7 < 23456$ and $1 > 23456$.
The correspondences for $P = 8$ are indicated with bold types in the tabulars of Figures~\ref{lsevenone}-\ref{lseventwo}, the other cases are easily deduced from these by the action of $D_8$.
We say that a list is {\em maximal\/} (or {\em minimal\/}) for some $F \in 1, \dots, 8$ if
$F$ lies outside (or inside) of all the conics determined by five of the other seven points.
Say $F = 8$, then for each possible $\hat 8 = n\pm$, there is one maximal and one minimal
list. The maximal list $\max(\hat 8 =  n\pm)$ is obtained from the given diagram $\hat 8 = n\pm$ taking  for all $\hat m$ with $m = 1, \dots 7$ the labels above the first column of horizontal arrows. 
The minimal list $\min(\hat 8 =  n\pm)$ is obtained from the given diagram $\hat 8 = n\pm$ taking  for all $\hat m$ with $m = 1, \dots 7$ the labels above the last column of horizontal arrows. These notations are consistant with $D_8$.
The extremal lists are easily realizable as follows: start from a configuration of eight points lying on the same conic $C_2$. Move first $8$, either to the interior or to the exterior of $C_2$. Then, move two further points $F$ and $G$, one to the interior and the other to the exterior of $C_2$, so as to obtain the desired list $\hat 8$. 
The extremal lists are distributed in 12 orbits: two are both maximal and minimal,
five are only maximal, five are only minimal. Note that $\max(\hat 1 = n\pm)$ and $\max(\hat 1 = (10-n)\mp)$ are deduced one from the other by $15$. We denote the orbit of these two lists shortly by $(n\pm, (10-n)\mp)$.  The maximal orbits $(2+, 8-)$ and $(3+, 7-)$ are also minimal:
$\max(\hat 1 = 8-) = \min(\hat 2 = 1-)$ and $\max(\hat 1 = 3+) = \min(\hat 2 = 1+)$.  
The principal extremal lists are: $L_{2} = \min(\hat 7 = 5-)$, $L_3 = \min(\hat 7 = 5+)$, $L_5 =  \min(\hat 7 = 3-)$, $L_{32} = \max(\hat 1 = 8+)$, $L_{48} = \max(\hat 1 = 6-)$, 
$L_{56} = \max(\hat 1 = 6+)$,  $L_{64} = \max(\hat 1 = 2+) = \min(\hat 8 = 1+)$, 
$L_{65} = \max(\hat 1 = 7-) = \min(\hat 8 = 1-)$, $L_{66} = \max(\hat 1 = 7+)$, 
$L_{67} = \max(\hat 1 = 5-)$, $L_{71} =  \min(\hat 1 = 4+)$, and $L_{72} = \min(\hat 1 = 4-)$.

\subsection{Another encoding for the lists and their orbits}

For the lists $L(1, \dots, 8)$, we define the {\em distance\/} between consecutive points $A, B$ in the cyclic ordering. Let $(A, \dots H) = (1, \dots 8)$ up to some cyclic permutation and reversion.
The distance $A \to B$ is the number of conics passing through
five points among $C, D, E, F, G, H$ that separate $A$ from $B$, multiplied by $-1$ if $A$ is the outermost of the two points.  
To determine a distance $A \to B$, one needs only to look at the sublists $\hat A$ and $\hat B$ of $L(1, \dots 8)$, as explained in the following example.
Let us find out the distance $8 \to 1$ in the two cases: $\hat 1 = 8+$, $\hat 8 = 1+$ and $\hat 1 = 8+$, $\hat 8 = 6+$. One has: $7 < 23456$, letting $8$ (resp. $1$) cross the sequence of six conics determined by $2, \dots 7$ yields the chain $\hat 1+(8)$ (resp. $\hat 8+(1)$), see Figure~\ref{chains}. For $\hat 1 = 8+, \hat 8 = 1+$, both points $8$ and $1$ lie outside of all six conics, $8 \to 1 = 0$, for $\hat 1 = 8+, \hat 8 = 6+$, one has $8 > 23456$, and $23467 < 1 < 23457$, thus $8 \to 1 = -2$. 

Note that the distance $A \to B$ is equal to $0$ if and only if: $\hat A =  B\epsilon$, $\hat B = A\epsilon$ or $\hat A = \hat B = N\epsilon$, with $N \in \{ C, \dots H  \}$.

\begin{figure}
\begin{displaymath}
\xymatrix{
8+ \ar[r]^{23456} &  6- \ar[r]^{23457} & 6+  \ar[r]^{23467} & 4- \ar[r]^{23567} & 4+ \ar[r]^{24567} & 2- 
 \ar[r]^{34567} & 2+ \\
1+ \ar[r]^{23456} &  6- \ar[r]^{23457} & 6+  \ar[r]^{23467} & 4- \ar[r]^{23567} & 4+ \ar[r]^{24567} & 2- 
 \ar[r]^{34567} & 2+ \\
}
\end{displaymath}
\caption{\label{chains}  The chains $\hat 1+(8)$ and $\hat 8+(1)$}
\end{figure}





The tabulars of Figures \ref{distone}-\ref{distwo} give the values of these invariants for all of the 47 successive principal lists $L_n$, and their sum $\sigma$.  
Note that the $47$ octuples $(1\to 2, \dots, 8 \to 1)$ listed in Figures~\ref{distone}-\ref{distwo}  are all different. Furthermore, their orbits under the action of $D_8$ are also all distinct. Observing the tabulars allows to derive the following:

\begin{proposition}
Any list $L(1, \dots, 8)$ is determined by the octuple $(1 \to 2, 2 \to 3, \dots, 8 \to 1)$.
Each orbit is determined by an octuple of integer numbers (ranging between $-6$ and $+6$), defined up to cyclic permutation, and reversion with change of all signs.
The absolute value $\vert \sigma \vert$ is an invariant of the orbits, that takes all of the even values between $0$ and $8$.  
The numbers of orbits realizing $\vert \sigma \vert = 0, 2, 4, 6, 8$ are respectively:  
$14$, $12$, $16$, $4$ and $1$.  
\end{proposition}

We have thus a new encoding for the lists and the orbits.

\begin{figure}
\begin{tabular}{ c c c c c c c c c }
$\hat 1$ & ${\bf 8+}$ & $8+$ & $8+$ & $8+$ & $8+$ & $8+$ & $8+$ & $8+$\\
$\hat 2$ & ${\bf 8+}$ & $8+$ & $8+$ & $8+$ & $8+$ & $8+$ & $8+$ & $8+$\\
$\hat 3$ & ${\bf 8+}$ & $8+$ & $8+$ & $8+$ & $8+$ & $8+$ & $8+$ & $8+$\\
$\hat 4$ & ${\bf 8+}$ & $8+$ & $8+$ & $8+$ & $6-$ & $6-$ & $6+$ & $3-$\\
$\hat 5$ & ${\bf 8+}$ & $8+$ & $6-$ & $6+$ & $6-$ & $6+$ & $6+$ & $3-$\\
$\hat 6$ & ${\bf 8+}$ & $5-$ & $5-$ & $5+$ & $5-$ & $5+$ & $3-$ & $3-$\\
$\hat 7$ & ${\bf 8+}$ & $5-$ & $5+$ & $5+$ & $3-$ & $3-$ & $3-$ & $3-$\\
 & & & & & & & &\\
$\hat 1$ & ${\bf 8+}$ & $8+$ & $8+$ & $8+$ & $8+$ & $8+$ & $8+$ & ${\bf 8+}$\\
$\hat 2$ & ${\bf 8+}$ & $8+$ & $8+$ & $8+$ & $8+$ & $8+$ & $8+$ & ${\bf 8+}$\\
$\hat 3$ & ${\bf 6-}$ & $6-$ & $6-$ & $6-$ & $6+$ & $6+$ & $4-$ & ${\bf 4+}$\\
$\hat 4$ & ${\bf 6-}$ & $6-$ & $6+$ & $3-$ & $6+$ & $3-$ & $3-$ & ${\bf 3+}$\\
$\hat 5$ & ${\bf 6-}$ & $6+$ & $6+$ & $3-$ & $6+$ & $3-$ & $3+$ & ${\bf 3+}$\\
$\hat 6$ & ${\bf 5-}$ & $5+$ & $3-$ & $3-$ & $3+$ & $3+$ & $3+$ & ${\bf 3+}$\\
$\hat 7$ & ${\bf 3+}$ & $3+$ & $3+$ & $3+$ & $3+$ & $3+$ & $3+$ & ${\bf 3+}$\\
 & & & & & & & &\\
$\hat 1$ & ${\bf 8+}$ & $8+$ & $8+$ & $8+$ & $8+$ & $8+$ & $8+$ & ${\bf 8+}$\\
$\hat 2$ & ${\bf 6-}$ & $6-$ & $6-$ & $6-$ & $6-$ & $6-$ & $6-$ & ${\bf 6-}$\\
$\hat 3$ & ${\bf 6-}$ & $6-$ & $6-$ & $6-$ & $6+$ & $6+$ & $4-$ & ${\bf 4+}$\\
$\hat 4$ & ${\bf 6-}$ & $6-$ & $6+$ & $3-$ & $6+$ & $3-$ & $3-$ & ${\bf 3+}$\\
$\hat 5$ & ${\bf 6-}$ & $6+$ & $6+$ & $3-$ & $6+$ & $3-$ & $3+$ & ${\bf 3+}$\\
$\hat 6$ & ${\bf 5-}$ & $5+$ & $3-$ & $3-$ & $3+$ & $3+$ & $3+$ & ${\bf 3+}$\\
$\hat 7$ & ${\bf 1-}$ & $1-$ & $1-$ & $1-$ & $1-$ & $1-$ & $1-$ & ${\bf 1-}$\\
 & & & & & & & &\\
$\hat 1$ & $8+$ & $8+$ & ${\bf 8+}$ & ${\bf 8+}$ & ${\bf 8+}$ & ${\bf 8+}$ & ${\bf 8+}$ & $8+$\\
$\hat 2$ & $6+$ & $6+$ & ${\bf 6+}$ & ${\bf 6+}$ & ${\bf 4-}$ & ${\bf 4-}$ & ${\bf 4+}$ & $1-$\\
$\hat 3$ & $6+$ & $6+$ & ${\bf 4-}$ & ${\bf 4+}$ & ${\bf 4-}$ & ${\bf 4+}$ & ${\bf 4+}$ & $1-$\\
$\hat 4$ & $6+$ & $3-$ & ${\bf 3-}$ & ${\bf 3+}$ & ${\bf 3-}$ & ${\bf 3+}$ & ${\bf 1-}$ & $1-$\\
$\hat 5$ & $6+$ & $3-$ & ${\bf 3+}$ & ${\bf 3+}$ & ${\bf 1-}$ & ${\bf 1-}$ & ${\bf 1-}$ & $1-$\\
$\hat 6$ & $1-$ & $1-$ & ${\bf 1-}$ & ${\bf 1-}$ & ${\bf 1-}$ & ${\bf 1-}$ & ${\bf 1-}$ & $1-$\\
$\hat 7$ & $1-$ & $1-$ & ${\bf 1-}$ & ${\bf 1-}$ & ${\bf 1-}$ & ${\bf 1-}$ & ${\bf 1-}$ & $1-$\\
\end{tabular}
\caption{\label{lone}  $\hat 8 = 1+$, lists $L_1, \dots, L_{32}$} 
\end{figure}

\begin{figure}
\begin{tabular}{ c c c c c c c c c }
$\hat 1$ & ${\bf 6-}$ & $6-$ & $6-$ & $6-$ & $6-$ & $6-$ & ${\bf 6-}$ & ${\bf 6-}$\\
$\hat 2$ & ${\bf 6-}$ & $6-$ & $6-$ & $6-$ & $6-$ & $6-$ & ${\bf 6-}$ & ${\bf 6-}$\\
$\hat 3$ & ${\bf 6-}$ & $6-$ & $6-$ & $6-$ & $6+$ & $6+$ & ${\bf 4-}$ & ${\bf 4+}$\\
$\hat 4$ & ${\bf 6-}$ & $6-$ & $6+$ & $3-$ & $6+$ & $3-$ & ${\bf 3-}$ & ${\bf 3+}$\\
$\hat 5$ & ${\bf 6-}$ & $6+$ & $6+$ & $3-$ & $6+$ & $3-$ & ${\bf 3+}$ & ${\bf 3+}$\\
$\hat 6$ & ${\bf 5-}$ & $5+$ & $3-$ & $3-$ & $3+$ & $3+$ & ${\bf 3+}$ & ${\bf 3+}$\\
$\hat 7$ & ${\bf 1+}$ & $1+$ & $1+$ & $1+$ & $1+$ & $1+$ & ${\bf 1+}$ & ${\bf 1+}$\\
 & & & & & & & &\\
$\hat 1$ & $6-$ & ${\bf 6-}$ & ${\bf 6-}$ & ${\bf 6-}$ & ${\bf 6-}$ & ${\bf 6-}$ & ${\bf 6-}$ & $6-$\\
$\hat 2$ & $6+$ & ${\bf 6+}$ & ${\bf 6+}$ & ${\bf 6+}$ & ${\bf 4-}$ & ${\bf 4-}$ & ${\bf 4+}$ & $1-$\\
$\hat 3$ & $6+$ & ${\bf 6+}$ & ${\bf 4-}$ & ${\bf 4+}$ & ${\bf 4-}$ & ${\bf 4+}$ & ${\bf 4+}$ & $1-$\\
$\hat 4$ & $6+$ & ${\bf 3-}$ & ${\bf 3-}$ & ${\bf 3+}$ & ${\bf 3-}$ & ${\bf 3+}$ & ${\bf 1-}$ & $1-$\\
$\hat 5$ & $6+$ & ${\bf 3-}$ & ${\bf 3+}$ & ${\bf 3+}$ & ${\bf 1-}$ & ${\bf 1-}$ & ${\bf 1-}$ & $1-$\\
$\hat 6$ & $1-$ & ${\bf 1-}$ & ${\bf 1-}$ & ${\bf 1-}$ & ${\bf 1-}$ & ${\bf 1-}$ & ${\bf 1-}$ & $1-$\\
$\hat 7$ & $1+$ & ${\bf 1+}$ & ${\bf 1+}$ & ${\bf 1+}$ & ${\bf 1+}$ & ${\bf 1+}$ & ${\bf 1+}$ & $1+$\\
 & & & & & & & &\\
$\hat 1$ & $6+$ & ${\bf 6+}$ & ${\bf 6+}$ & ${\bf 6+}$ & ${\bf 6+}$ & ${\bf 6+}$ & ${\bf 6+}$ & $6+$\\
$\hat 2$ & $6+$ & ${\bf 6+}$ & ${\bf 6+}$ & ${\bf 6+}$ & ${\bf 4-}$ & ${\bf 4-}$ & ${\bf 4+}$ & $1-$\\
$\hat 3$ & $6+$ & ${\bf 6+}$ & ${\bf 4-}$ & ${\bf 4+}$ & ${\bf 4-}$ & ${\bf 4+}$ & ${\bf 4+}$ & $1-$\\
$\hat 4$ & $6+$ & ${\bf 3-}$ & ${\bf 3-}$ & ${\bf 3+}$ & ${\bf 3-}$ & ${\bf 3+}$ & ${\bf 1-}$ & $1-$\\
$\hat 5$ & $6+$ & ${\bf 3-}$ & ${\bf 3+}$ & ${\bf 3+}$ & ${\bf 1-}$ & ${\bf 1-}$ & ${\bf 1-}$ & $1-$\\
$\hat 6$ & $1+$ & ${\bf 1+}$ & ${\bf 1+}$ & ${\bf 1+}$ & ${\bf 1+}$ & ${\bf 1+}$ & ${\bf 1+}$ & $1+$\\
$\hat 7$ & $1+$ & ${\bf 1+}$ & ${\bf 1+}$ & ${\bf 1+}$ & ${\bf 1+}$ & ${\bf 1+}$ & ${\bf 1+}$ & $1+$\\
 & & & & & & & &\\
$\hat 1$ & ${\bf 4-}$ & ${\bf 4-}$ & ${\bf 4-}$ & ${\bf 4-}$ & ${\bf 4+}$ & ${\bf 4+}$ & ${\bf 2-}$ & $2+$\\
$\hat 2$ & ${\bf 4-}$ & ${\bf 4-}$ & ${\bf 4+}$ & ${\bf 1-}$ & ${\bf 4+}$ & ${\bf 1-}$ & ${\bf 1-}$ & $1+$\\
$\hat 3$ & ${\bf 4-}$ & ${\bf 4+}$ & ${\bf 4+}$ & ${\bf 1-}$ & ${\bf 4+}$ & ${\bf 1-}$ & ${\bf 1+}$ & $1+$\\
$\hat 4$ & ${\bf 3-}$ & ${\bf 3+}$ & ${\bf 1-}$ & ${\bf 1-}$ & ${\bf 1+}$ & ${\bf 1+}$ & ${\bf 1+}$ & $1+$\\
$\hat 5$ & ${\bf 1+}$ & ${\bf 1+}$ & ${\bf 1+}$ & ${\bf 1+}$ & ${\bf 1+}$ & ${\bf 1+}$ & ${\bf 1+}$ & $1+$\\
$\hat 6$ & ${\bf 1+}$ & ${\bf 1+}$ & ${\bf 1+}$ & ${\bf 1+}$ & ${\bf 1+}$ & ${\bf 1+}$ & ${\bf 1+}$ & $1+$\\
$\hat 7$ & ${\bf 1+}$ & ${\bf 1+}$ & ${\bf 1+}$ & ${\bf 1+}$ & ${\bf 1+}$ & ${\bf 1+}$ & ${\bf 1+}$ & $1+$\\
\end{tabular}
\caption{\label{ltwo} $\hat 8 = 1+$, lists $L_{33}, \dots L_{64}$} 
\end{figure}

\begin{figure}
\begin{tabular}{ c c c c c c c }
$\hat 1$ & $7-$ & $7+$ & $5-$ & ${\bf 5+}$ & ${\bf 3-}$ & ${\bf 3+}$\\
$\hat 2$ & $1+$ & $1+$ & $1+$ & ${\bf 1+}$ & ${\bf 1+}$ & ${\bf 1+}$\\
$\hat 3$ & $1+$ & $1+$ & $1+$ & ${\bf 1+}$ & ${\bf 1+}$ & ${\bf 1-}$\\
$\hat 4$ & $1+$ & $1+$ & $1+$ & ${\bf 1+}$ & ${\bf 1-}$ & ${\bf 1-}$\\
$\hat 5$ & $1+$ & $1+$ & $1+$ & ${\bf 1-}$ & ${\bf 1-}$ & ${\bf 1-}$\\
$\hat 6$ & $1+$ & $1+$ & $1-$ & ${\bf 1-}$ & ${\bf 1-}$ & ${\bf 1-}$\\
$\hat 7$ & $1+$ & $1-$ & $1-$ & ${\bf 1-}$ & ${\bf 1-}$ & ${\bf 1-}$\\
\end{tabular}
\caption{\label{lthree}  $\hat 8 = 1-$, lists $L_{65}, \dots, L_{70}$} 
\end{figure}

\begin{figure}
\begin{tabular}{ c c c c c }
$\hat 1$ & $4+$ & $4-$ & ${\bf 6+}$ & ${\bf 6-}$\\
$\hat 2$ & $8-$ & $8-$ & ${\bf 8-}$ & ${\bf 8-}$\\
$\hat 3$ & $8-$ & $8-$ & ${\bf 8-}$ & ${\bf 8-}$\\
$\hat 4$ & $2+$ & $8-$ & ${\bf 8-}$ & ${\bf 8-}$\\
$\hat 5$ & $2+$ & $2+$ & ${\bf 8-}$ & ${\bf 8-}$\\
$\hat 6$ & $2+$ & $2+$ & ${\bf 2+}$ & ${\bf 8-}$\\
$\hat 7$ & $2+$ & $2+$ & ${\bf 2+}$ & ${\bf 2+}$\\
\end{tabular}
\caption{\label{lfour} $\hat 8 = 2+$, lists $L_{71}, \dots, L_{74}$} 
\end{figure}

\begin{figure}
\begin{tabular}{ c c c c c c c c }
$\hat 1$ & $4+$ & ${\bf 4+}$ & ${\bf 6-}$ & $4+$ & ${\bf 4+}$ &  $4-$ & ${\bf 6+}$\\
$\hat 2$ & $8+$ & ${\bf 4-}$ & ${\bf 8+}$ & $6-$ & ${\bf 6+}$ &  $8+$ & ${\bf 8+}$\\
$\hat 3$ & $8-$ & ${\bf 8-}$ & ${\bf 8-}$ & $8-$ & ${\bf 8-}$ &  $8-$ & ${\bf 8-}$\\
$\hat 4$ & $2+$ & ${\bf 2+}$ & ${\bf 8-}$ & $2+$ & ${\bf 2+}$ &  $8-$ & ${\bf 8-}$\\
$\hat 5$ & $2+$ & ${\bf 2-}$ & ${\bf 8-}$ & $2+$ & ${\bf 2+}$ &  $2+$ & ${\bf 8-}$\\
$\hat 6$ & $2+$ & ${\bf 2-}$ & ${\bf 8-}$ & $2+$ & ${\bf 2-}$ &  $2+$ & ${\bf 2+}$\\
$\hat 7$ & $2+$ & ${\bf 2-}$ & ${\bf 2+}$ & $2-$ & ${\bf 2-}$ &  $2+$ & ${\bf 2+}$\\
 & & & & &  & & \\
$\hat 1$ &  $6+$ & $6+$ & $4-$ & ${\bf 6-}$ & $4-$ & $4-$ & \\
$\hat 2$ &  $6-$ & $6+$ & $4-$ & ${\bf 6-}$ & $6-$ & $6+$ & \\
$\hat 3$ &  $8-$ & $8-$ & $8-$ & ${\bf 8-}$ & $8-$ & $8-$ & \\
$\hat 4$ &  $8-$ & $8-$ & $8-$ & ${\bf 8-}$ & $8-$ & $8-$ & \\
$\hat 5$ &  $8-$ & $8-$ & $2-$ & ${\bf 8-}$ & $2+$ & $2+$ & \\
$\hat 6$ &  $2+$ & $2-$ & $2-$ & ${\bf 8-}$ & $2+$ & $2-$ & \\
$\hat 7$ &  $2-$ & $2-$ & $2-$ & ${\bf 2-}$ & $2-$ & $2-$ & \\
\end{tabular}
\caption{\label{lfive}  $\hat 8 = 2-$, lists $L_{75}, \dots, L_{87}$}
\end{figure}

\begin{figure}
\begin{tabular}{ c c c c c c c c c }
$\hat 1$ & $5+$ & ${\bf 5+}$ & ${\bf 5-}$ & ${\bf 5+}$ & $5-$ & ${\bf 5-}$ & 
${\bf 5+}$ & ${\bf 5-}$\\
$\hat 2$ & $7+$ & ${\bf 5-}$ & ${\bf 7+}$ & ${\bf 5-}$ & $7+$ & ${\bf 5-}$ & ${\bf 7+}$ & ${\bf 5-}$\\
$\hat 3$ & $7-$ & ${\bf 7+}$ & ${\bf 7-}$ & ${\bf 7-}$ & $7+$ & ${\bf 7-}$ & 
${\bf 7+}$ & ${\bf 7+}$\\
$\hat 4$ & $1+$ & ${\bf 1+}$ & ${\bf 1+}$ & ${\bf 1+}$ & $1+$ & ${\bf 1+}$ & ${\bf 1+}$ & ${\bf 1+}$\\
$\hat 5$ & $1-$ & ${\bf 1-}$ & ${\bf 1+}$ & ${\bf 1-}$ & $1+$ & ${\bf 1+}$ & 
${\bf 1-}$ & ${\bf 1+}$\\
$\hat 6$ & $1-$ & ${\bf 3+}$ & ${\bf 1-}$ & ${\bf 3+}$ & $1-$ & ${\bf 3+}$ & 
${\bf 1-}$ & ${\bf 3+}$\\
$\hat 7$ & $3+$ & ${\bf 3-}$ & ${\bf 3+}$ & ${\bf 3+}$ & $3-$ & ${\bf 3+}$ & ${\bf 3-}$ & ${\bf 3-}$\\
\end{tabular}
\caption{\label{lsix}  $\hat 8 = 3-$, lists $L_{88}, \dots, L_{95}$}
\end{figure}

\begin{figure}
\begin{tabular}{|c |c|c|c|c|c|c|c|c|c|}
\hline
$n$ & $1 \to 2$ & $2 \to 3$ & $3 \to 4$ & $4 \to 5$ & $5 \to 6$ & $6 \to 7$ & $7 \to 8$ & $8 \to 1$ & $\sigma$ \\
\hline
$2$ & $0$ & $0$ & $0$ & $0$ & $-1$ & $0$ & $5$ & $0$ & $4$ \\
$3$ & $0$ & $0$ & $0$ & $1$ & $0$ & $-1$ & $4$ & $0$ & $4$ \\
$4$ & $0$ & $0$ & $0$ & $2$ & $0$ & $0$ & $4$ & $0$ & $6$ \\
$5$ & $0$ & $0$ & $-1$ & $0$ & $0$ & $-2$ & $3$ & $0$ & $0$ \\
$6$ & $0$ & $0$ & $-1$ & $1$ & $0$ & $-1$ & $3$ & $0$ & $2$ \\
$7$ & $0$ & $0$ & $-2$ & $0$ & $1$ & $0$ & $3$ & $0$ & $2$ \\
$8$ & $0$ & $0$ & $-3$ & $0$ & $0$ & $0$ & $3$ & $0$ & $0$ \\
$10$ & $0$ & $1$ & $0$ & $1$ & $0$ & $-2$ & $2$ & $0$ & $2$ \\
$11$ & $0$ & $1$ & $-1$ & $0$ & $1$ & $-1$ & $2$ & $0$ & $2$ \\
$12$ & $0$ & $1$ & $-2$ & $0$ & $0$ & $-1$ & $2$ & $0$ & $0$ \\
$13$ & $0$ & $2$ & $0$ & $0$ & $2$ & $0$ & $2$ & $0$ & $6$ \\
$14$ & $0$ & $2$ & $-1$ & $0$ & $1$ & $0$ & $2$ & $0$ & $4$ \\
$15$ & $0$ & $3$ & $0$ & $-1$ & $0$ & $0$ & $2$ & $0$ & $4$ \\
$18$ & $-1$ & $0$ & $0$ & $1$ & $0$ & $-3$ & $1$ & $0$ & $-2$ \\
$19$ & $-1$ & $0$ & $-1$ & $0$ & $1$ & $-2$ & $1$ & $0$ & $-2$ \\
$20$ & $-1$ & $0$ & $-2$ & $0$ & $0$ & $-2$ & $1$ & $0$ & $-4$ \\
$21$ & $-1$ & $1$ & $0$ & $0$ & $2$ & $-1$ & $1$ & $0$ & $2$ \\
$22$ & $-1$ & $1$ & $-1$ & $0$ & $1$ & $-1$ & $1$ & $0$ & $0$ \\
$23$ & $-1$ & $2$ & $0$ & $-1$ & $0$ & $-1$ & $1$ & $0$ & $0$ \\
$25$ & $-2$ & $0$ & $0$ & $0$ & $3$ & $0$ & $1$ & $0$ & $2$ \\
$26$ & $-2$ & $0$ & $-1$ & $0$ & $2$ & $0$ & $1$ & $0$ & $0$ \\
$32$ & $-5$ & $0$ & $0$ & $0$ & $0$ & $0$ & $1$ & $0$ & $-4$ \\
$34$ & $0$ & $0$ & $0$ & $1$ & $0$ & $-4$ & $0$ & $1$ & $-2$ \\
$35$ & $0$ & $0$ & $-1$ & $0$ & $1$ & $-3$ & $0$ & $1$ & $-2$ \\
$36$ & $0$ & $0$ & $-2$ & $0$ & $0$ & $-3$ & $0$ & $1$ & $-4$ \\
$37$ & $0$ & $1$ & $0$ & $0$ & $2$ & $-2$ & $0$ & $1$ & $2$ \\
$38$ & $0$ & $1$ & $-1$ & $0$ & $1$ & $-2$ & $0$ & $1$ & $0$ \\
$41$ & $-1$ & $0$ & $0$ & $0$ & $3$ & $-1$ & $0$ & $1$ & $2$ \\
$48$ & $-4$ & $0$ & $0$ & $0$ & $0$ & $-1$ & $0$ & $1$ & $-4$ \\
$49$ & $0$ & $0$ & $0$ & $0$ & $4$ & $0$ & $0$ & $2$ & $6$ \\
$56$ & $-3$ & $0$ & $0$ & $0$ & $1$ & $0$ & $0$ & $2$ & $0$ \\
$64$ & $0$ &  $0$ & $0$ & $0$ & $0$ & $0$ & $0$ & $6$ & $6$ \\
\hline
\end{tabular}
\caption{\label{distone} Distances $A \to B$ for the principal lists with $\hat 8 = 1+$} 
\end{figure}

\begin{figure}
\begin{tabular}{|c|c|c|c|c|c|c|c|c|c|}
\hline
$n$ & $1 \to 2$ & $2 \to 3$ & $3 \to 4$ & $4 \to 5$ & $5 \to 6$ & $6 \to 7$ & $7 \to 8$ & $8 \to 1$ & $\sigma$ \\
\hline
$65$ & $-1$ & $0$ & $0$ & $0$ & $0$ & $0$ & $-1$ & $6$ & $4$ \\
$66$ & $-2$ & $0$ & $0$ & $0$ & $0$ & $1$ & $0$ & $5$ & $4$ \\
$67$ & $-3$ & $0$ & $0$ & $0$ & $-1$ & $0$ & $0$ & $4$ & $0$ \\
$71$ & $5$ & $0$ & $1$ & $0$ & $0$ & $0$ & $0$ & $-2$ & $4$ \\
$72$ & $4$ & $0$ & $0$ & $-1$ & $0$ & $0$ & $0$ & $-3$ & $0$ \\ 
$75$ & $4$ & $-1$ & $1$ & $0$ & $0$ & $0$ & $1$ & $-1$ & $4$ \\
$78$ & $3$ & $-2$ & $1$ & $0$ & $0$ & $-1$ & $0$ & $-1$ & $0$ \\
$80$ & $3$ & $-1$ & $0$ & $-1$ & $0$ & $0$ & $1$ & $-2$ & $0$ \\
$82$ & $1$ & $-2$ & $0$ & $0$ & $1$ & $-1$ & $0$ & $-3$ & $-4$ \\
$83$ & $0$ & $-3$ & $0$ & $0$ & $2$ & $0$ & $0$ & $-3$ & $-4$ \\
$84$ & $0$ & $-4$ & $0$ & $-2$ & $0$ & $0$ & $0$ & $-2$ & $-8$ \\
$86$ & $2$ & $-2$ & $0$ & $-1$ & $0$ & $-1$ & $0$ & $-2$ & $-4$ \\
$87$ & $1$ & $-3$ & $0$ & $-1$ & $1$ & $0$ & $0$ & $-2$ & $-4$ \\
$88$ & $-2$ & $1$ & $-1$ & $1$ & $0$ & $1$ & $-1$ & $1$ & $0$ \\
$92$ & $-1$ & $0$ & $-2$ & $0$ & $-1$ & $2$ & $0$ & $2$ & $0$ \\
\hline
\end{tabular}
\caption{\label{distwo}  Distances $A \to B$ for the principal lists with $\hat 8 \not= 1+$} 
\end{figure}

\subsection{Isotopies of octuples of points}

Consider the space $(\mathbb{R}P^2)^8$ stratified by the lines and conics.

\begin{proposition}
The octuples of points realizing a given list $L(1, \dots 8)$ lie all in the same chamber of $(\mathbb{R}P^2)^8$. This chamber is adjacent to the deep stratum, formed by the configurations of eight coconic points. 
\end{proposition}

\begin{description}
\item[Proof:]
Let $1, \dots 8$ realize some list $L$.
Denote by $h_t$ a homothety centered at any point $P$ different from $1, \dots 8$, with rate $t$.
The homotheties $h_t, t \in [ 1, \infty ]$ give rise to an isotopy of $1, \dots 8$.
During this isotopy, the list $L$ is preserved, except for the end-point when the mobile points $1, \dots 8$ become all aligned. They are also on some reducible conic. The configuration is thus on the deep stratum. Move in this stratum so that the conic becomes non-reducible. With a slight perturbation, we can make the isotopy generic, in other words, it is completely contained in the same chamber except for the end-point.
Consider now two configurations realizing the same list. Perform for each of them a generic isotopy till all of the eight points lie on the same conic. Up to an affine transformation mapping one conic onto the other and an isotopy along this conic, we may assume that both isotopies have the same end-point.
So we have a path connecting the two configurations and having one single non-generic point. Perturb the path in a neighbourhood of this point. In this neighbourhood, all of the walls intersect pairwise transversally. There are two possible perturbations, one crosses all of the walls twice, the other crosses no wall at all. $\Box$
\end{description}

\subsection{Elementary changes and inductive constructions}

\begin{figure}[htbp]
\centering
\includegraphics{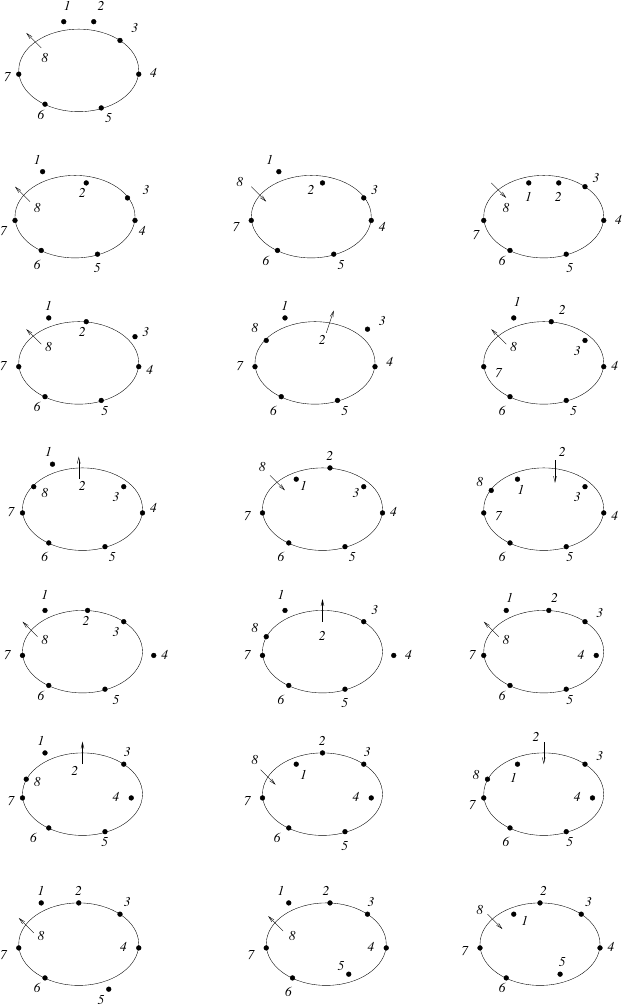}
\caption{\label{elem}  Elementary changes}
\end{figure}

In this section, we finish the proof of Theorem~2.
Let us call {\em elementary change\/} the change induced on a list $L = L(1, \dots 8)$ letting one point among $1, \dots, 8$ cross a conic determined by five others, in some direction. Up to the action of $D_8$, there are 19 elementary changes, see Figures~\ref{elem}-\ref{repres}. 
A pair $(\hat N = Q\pm, \hat M = R\pm)$ involved in some elementary change is called {\em elementary pair\/}. There are in total $224$ elementary changes and correspondingly $224$ elementary pairs, see Figures~\ref{firstel}-\ref{secondel}. 

\begin{figure}
\begin{tabular}{c c c c c}
$\hat 1: 2+ \to 2-$ & & & & \\
$\hat 2: 1+ \to 1-$ & & & & \\
 & & \\
$\hat 1: 3+ \to 8-$ & & $\hat 1: 8- \to 3+$ & & $\hat 1: 8- \to 3+$\\
$\hat 2: 1+ \to 1-$ & & $\hat 2: 1- \to 1+$ & & $\hat 2: 8- \to 3+$\\
 & & \\
$\hat 1: 3- \to 3+$ & & $\hat 1: 3+ \to 3-$ & & $\hat 1: 2- \to 4+$\\
$\hat 3: 1+ \to 1-$ & & $\hat 3: 1- \to 1+$ & & $\hat 3: 1+ \to 1-$\\
 & & \\
$\hat 1: 4+ \to 2-$ & & $\hat 1: 4+ \to 2-$ & & $\hat 1: 2- \to 4+$\\
$\hat 3: 1- \to 1+$ & & $\hat 3: 8- \to 2+$ & & $\hat 3: 2+ \to 8-$\\
 & & \\
$\hat 1: 4+ \to 4-$ & & $\hat 1: 4- \to 4+$ & & $\hat 1: 5+ \to 3-$\\
$\hat 4: 1+ \to 1-$ & & $\hat 4: 1- \to 1+$ & & $\hat 4: 1+ \to 1-$\\
 & & \\
$\hat 1: 3- \to 5+$ & & $\hat 1: 3- \to 5+$ & & $\hat 1: 5+ \to 3-$\\
$\hat 4: 1- \to 1+$ & & $\hat 4: 8- \to 2+$ & & $\hat 4: 2+ \to 8-$\\
 & & \\
$\hat 1: 5- \to 5+$ & & $\hat 1: 4- \to 6+$ & & $\hat 1: 6+ \to 4-$\\
$\hat 5: 1+ \to 1-$ & & $\hat 5: 1+ \to 1-$ & & $\hat 5: 8- \to 2+$\\
\end{tabular}
\caption{\label{repres} Representants of the 19 orbits of elementary changes}
\end{figure}

\begin{proposition}
Let $1, \dots, 8$ be eight points in convex position, realizing a list L. Let $A, B$ be two of these points, consecutive for the cyclic ordering.
\begin{enumerate}
\item
If $A \to B = 0$, one may move the set of points, shifting $A$ towards $B$, till $A = B$, without degeneration of the list $L(1, \dots 8)$ inbetween.
\item
If $A \to B = \pm n$, with $n > 0$, the points $A$ and $B$ are separated from each other by $n$ of the six conics determined by $C, D, E, F, G, H$. One may move the set of points, shifting $A$ towards $B$, till $A = B$, in such a way that the list $L(1, \dots, 8)$ undergoes exactly $n$ degenerations, corresponding to the $n$ conics.
\end{enumerate}
\end{proposition}

\begin{description}
\item[Proof:]
For any pair of consecutive points $P, Q$ among $1, \dots 8$, denote by $[PQ]$ the affine segment $PQ$ belonging to the convex hull of the eight points. Move $A$ towards $B$ as prescribed, leaving the other points fixed. The point $A$ moves inside of a triangle $BHR$ supported by the lines $(BC)$, $(BH)$, $(GH)$, where $R = (BC) \cap (GH)$. In most cases, on may find a path connecting the starting point $A$ to $B$ inside of the triangle, that doesn't cut any one of the $15$ conics passing through $B$ and four other points among $C, \dots H$.
The exception is the following: let $C_2$ be one of the conics $BDEFG$, $BDEFH$, $BDEGH$, $BDFGH$ or $BEFGH$. If $A, C$ lie both outside of $C_2$, and the second intersection of $C_2$ with the line $BC$ is on the edge $BR$ of the triangle, then any path $8 \to 1$ must cross $C_2$.  
Assume now that we may move all of the points at the same time. We perform an isotopy in $(\mathbb{R}P^2)^8$ inside of the chamber of $L$, towards the deep stratum. When the octuple of points gets close enough to this stratum, the exceptional case can no longer occur: let $C_2$ be any one of the five conics, and assume that $A, C$ lie both outside of $C_2$, then the second intersection of $C_2$ with the line $(BC)$ lies on the segment $[BC]$. $\Box$. 
\end{description}

\begin{proposition}
The elementary changes are always realizable: for each elementary pair appearing in 
an existing list $L$, one may perform the corresponding elementary change to  realize a new list.
\end{proposition}

\begin{description}
\item[Proof:]
Recall that a statement of the form {\em $6$ crosses $12345$ from the inside to the outside\/} is equivalent to {\em $5$ crosses $12346$ from the outside to the inside\/}.
All of the $19$ elementary changes but one may be thus interpreted as motions of a point $A$ towards a consecutive point $B$ such that $A \to B \not= 0$, till $A$ crosses the first conic separating it from $B$.
The only exception is the third change: $\hat 1: 8- \to 3+$, $\hat 2: 1- \to 1+$. So, for any elementary pair of $L$ that is not in the orbit of  ($\hat 1 = 8-$, $\hat 2 = 1-$), the corresponding elementary change is realizable. Let now $L$ have
the elementary pair ($\hat 1 = 8-$, $\hat 2 = 1-$).
Note that $\hat 1 = 8-$ belongs to the chain $\hat 1-(8)$, and $\hat 2 = 1-$ belongs to the chain $\hat 2+(8)$. Using Figure~\ref{values}, we deduce that $\hat 8 = 1-$; using  Figure~\ref{oneminus}, we see that $\hat N = 1-$ for $N = 3, \dots 7$. The list is $\max(\hat 1 = 8-)$.
The point $8$ lies inside of all the conics determined by $1$ and four points among $2, \dots 7$. Up to some isotopy towards the deep stratum, we may assume that the conic $34567$ has its second intersection point with the line $(17)$ on the affine segment $[17]$. Denote by $T$ the triangle supported by the lines $(67)$, $(71)$, $(12)$ containing $8$.
One may move $8$ towards $[17]$ in $T$, preserving the list, till $8$ reaches $34567$. $\Box$
\end{description} 

Start from a list that is already realized, we will say that we perform the change $(\hat N, \hat M)$ without further precision, as there is no ambiguity possible. 
Any principal list may be obtained from an extremal list by some sequence of elementary changes and actions
of elements of $D_8$. See Figure~\ref{induct} where $n$ stands for a list $L_n$  (as in Figures~\ref{distone}-\ref{distwo}) and each row is a new sequence whose first list was already realized.
These sequences are chosen so as to reach all of the principal lists with the least possible number of starting lists 
(note that some intermediate lists are also extremal). This finishes the proof of Theorem~2. $\Box$

\begin{figure}
\begin{tabular}{ c c c c c c c c c c c c c }
$1$ & $(\hat 6, \hat 7)$  & $2$ & $(\hat 5, \hat 7)$ & $3$ & $(\hat 5, \hat 6)$ & $4$ & $(\hat 4, \hat 7)$  & $6$ & $(\hat 4, \hat 6)$ & $7$ & $(\hat 4, \hat 5)$ & $8$  \\
$6$ & $(\hat 5, \hat 6)$ & $5$ & $(37)$ & $\bf{9}$ & & & & & & & & \\
$\bf{9}$ & $(\hat 5, \hat 6)$ & $10$ & $(\hat 4,  \hat 6)$ & $11$ & $(\hat 4, \hat 5)$ & $12$ & $(\hat 3, \hat 6)$ & $14$ & $(\hat 4, \hat 5)$ & $13$ & &  \\
$14$ & $(\hat 3, \hat 5)$ & $15$ & & & & & & & & & & \\
$3$ & $(37)$ & $\bf{17}$ & $(\hat 5, \hat 6)$ & $18$ & $(\hat 4, \hat 6)$ & $19$ & $(\hat 4, \hat 5)$ & 20 & $(\hat 3, \hat 6)$ & $22$ & $(\hat 3, \hat 5)$ & $23$\\
$22$ & $(\hat 4, \hat 5)$ & $21$ & $(\hat 2, \hat 6)$ & $25$ & $(\hat 4, \hat 5)$ & $26$ & & & & & & \\
$2$ & $(37)$ & $\bf{33}$ & $(\hat 5, \hat 6)$ & $34$ & $(\hat 4, \hat 6)$ & $35$ & $(\hat 4,  \hat 5)$ & $36$ & $(\hat 3, \hat 6)$ & $38$ & & \\
$35$ & $(\hat 3, \hat 6)$ & $37$ & $(\hat 2, \hat 6)$ & $41$ & $(\hat 1, \hat 6)$ & $49$ & & & & & & \\
$65$ & $(\hat 1, \hat 7)$ & $66$ & $(\hat 1, \hat 6)$ & $67$ & & & & & & & & \\
$71$ & $(\hat 1, \hat 4)$ & $72$ & & & & & & & & & & \\
$71$ & $(\hat 2, \hat 8)$ & $75$ & $(\hat 2, \hat 7)$ & $78$ & & & & & & & & \\
$75$ & $(\hat 1, \hat 4)$ & $80$ & $(\hat 2, \hat 7)$ & $86$ & $(\hat 2, \hat 6)$ & $87$ & $(\hat 2, \hat 5)$ & $84$ & & & & \\
$87$ & $(\hat 1, \hat 5)$ & $83$ & $(\hat 2, \hat 6)$ & $82$ & & & & & & & & \\
$87$ & $(\hat 3, \hat 8)$ & $(-1)$ & $88$ & $(15)$  & $\bf{90}$ & $(\hat 3, \hat 7)$ & $92$ & & & & & \\
 & & & or & $(\hat 1, \hat 5)$ & & & & & & & &\\
\end{tabular}
\caption{\label{induct} Inductive construction of the principal lists}
\end{figure}

\section{Link between lists and pencils}

\subsection{Nodal lists}

\begin{figure}[htbp]
\centering
\includegraphics{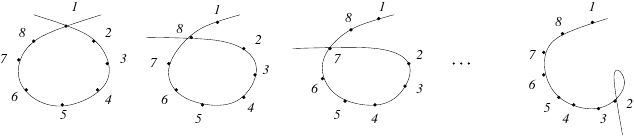}
\caption{\label{combi} The cubics $(1\pm, 1)_{nod}$, $(1-, 8)_{nod}, \dots (1-, 2)_{nod}$}
\end{figure}

A configuration of eight points lying on a nodal cubic $C_3$, one of them being the node, is not generic. Up to the action of $D_8$, there are eight possible combinatorial cubics $C_3$, see Figure~\ref{combi} where the successive
cubics are denoted by $(1\pm, 1)_{nod}$,  $(1-, k)_{nod}$ for $k = 8, \dots 2$.
Note that the encoding is consistant with the action of $D_8$.

\begin{definition}
A list $L(1, \dots, 8)= (\hat 1, \hat 2, \dots, \hat 8)$ is {\em nodal\/} if it is realizable by $1, \dots, 8$ on a nodal cubic with double point at one of these points.
\end{definition}

\begin{proposition}
Up to the action of the group $D_8$, 
$\natural L(1, \dots, 8)_{nodal}=4$. The nodal orbits are the maximal orbits $(8-, 2+)$, $(8+, 2-)$, $(6-, 4+)$, $(6+, 4-)$.  
Here are representants of each orbit along with the corresponding nodal cubics.

$L(1, \dots, 8)= (\hat 1, \hat 2, \dots, \hat 8)=$

$\max(\hat 1 = 8-)$, realizable with $(1\pm, 1)_{nod}$, $(1-, k)_{nod}, k = 2, \dots 8$,

$\max(\hat 1 = 8+)$, realizable with $(1\pm, 1)_{nod}$ and $(1-, 8)_{nod}$,

$\max(\hat 1 = 6-)$ and $\max(\hat 1 = 6+)$, both realizable with $(1\pm, 1)_{nod}$.
\end{proposition}

\begin{description}
\item[Proof:]
Let $1, \dots, 8$ lie in convex position on some nodal cubic $C_3$, one of 
these points being the node. Up to the action of $D_8$, we may assume that
$C_3$ is one of the cubics $(1\pm, 1)_{nod}$, $(1-, k)_{nod}, k = 2 \dots 8$. 
If $C_3 = (1\pm, 1)_{nod}$, one has $\hat i \in \{1\pm\}$ for all 
$i \in \{2, \dots, 8\}$. 
Using the two diagrams $\hat 8 =1\pm$ of Figures~\ref{oneplus}-\ref{oneminus} , one finds 14 possibilities for the list $L(1, \dots, 8)$, namely the maximal lists $\max(\hat 1 = n\pm)$, $n = 2, \dots, 8$.
One has $\hat 1 = 8+ \iff 7 < 23456$ and $8 > 23456$. One can choose the 
points $2, \dots, 8$ on the loop of $C_3$ so that this condition is 
achieved. One has $\hat 1=7-  \iff 8 < 23456$ and 
$7 > 23456$. By Bezout's theorem between $23456$ and $(1\pm, 1)_{nod}$, one 
cannot choose the points $1, \dots, 8$ on the loop of this cubic verifying 
this condition. Finishing this argument with the other possible values of 
$\hat 1$, one finds that $1, \dots, 8$ may be chosen on $(1\pm, 1)_{nod}$ so as to realize the eight lists $\max(\hat 1 = n\pm)$, with $n = 2, 4, 6, 8$.

Similarly, one proves that $1, \dots, 8$ on the cubic $C_3 = (1-, 8)_{nod}$ can realize 
exactly the first two lists $\max(\hat 1 = 8\pm)$; and points $1, \dots, 8$ on 
any cubic $C_3 = (1-, k)_{nod}, k \in \{2, \dots, 7\}$ must realize the first list $\max(\hat 1 = 8-)$. $\Box$
\end{description}

Given two points $P, Q$ among $1, \dots, 8$, we denote by $(\hat P, Q)$ the cubic of the list $\hat P$ with double point at $Q$.  

\begin{proposition}
The three conditions hereafter are equivalent:
\begin{enumerate}
\item
The list $L(1, \dots, 8)$ determines the (combinatorial) pencil 
$\mathcal{P}(1, \dots,8)$,
\item
The list $L(1, \dots, 8)$ is not nodal,
\item
$\forall G \in \{1, \dots, 8\}, \forall C_3$ cubic of $\hat G$, the 
position of $G$ with respect to $C_3$ is determined by $L(1, \dots, 8)$.
\end{enumerate}
\end{proposition}

\begin{figure}[htbp]
\centering 
\includegraphics{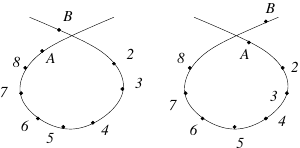}
\caption{\label{nodal}  $C^1_3$, $C^2_3$: $(A, B) = (1, 9)$, $C^3_3$, $C^4_3$: $(A, B) = (9, 1)$}
\end{figure}

\begin{description}
\item[Proof:]
$\neg 2 \Rightarrow \neg 1$: let $L(1, \dots, 8)$ be a nodal list. Up to the action of $D_8$, we may assume it is one of the four lists in Proposition~6. Any one of these lists is realizable with the eight points on a nodal cubic $C_3 = (1\pm, 1)_{nod}$. The points $1, \dots, 8$ give rise to a singular pencil
$\mathcal{P}_{sing}$ with $1 = 9$. Perturb the pencil moving $1$ away from the node onto the odd component or onto the loop of $C_3$ (leaving the other seven points fixed). In the generic pencil obtained, $C_3$ is replaced by a pair of distinguished cubics $C^1_3$, $C^2_3$ or $C^3_3$, $C^4_3$, see Figure~\ref{nodal}. Using Bezout's theorem, we see easily that these two pairs of
cubics cannot belong to the same pencil. 

$2 \Rightarrow 3$: according to Proposition~3, any two configurations of points realizing the same list may be connected by a path inside of their common chamber of $(\mathbb{R}P^2)^8$. As $L(1, \dots, 8)$ is not nodal, none of the points $G$ may cross any cubic of the list $\hat G$.

$3 \Rightarrow 1$: consider a configuration of points realizing a list $L(1, \dots 8)$ and giving rise to some pencil $\mathcal{P}$. Move the configuration preserving the list, the combinatorial pencil degenerates only if some base point $A$ among $1, \dots 8$ comes together with $9$. But this amounts to say that each point $G \not= A \in \{ 1, \dots, 8 \}$ comes onto the cubic $(\hat G, A)$. 
$\Box$
\end{description}

Propositions~4 and 7 imply immediately: 

\begin{proposition}
Let $1, \dots, 8$ be eight points in convex position, realizing a non-nodal list L. Let $A, B$ be two of these points, consecutive for the cyclic ordering.
We may move the set of points, shifting $A$ towards $B$, in such a way that
the combinatorial pencil does not degenerate till $A = B$ (if $A \to B = 0$)
or till $A$ reaches the first of the $n$ conics separating $A$ from $B$ (if $A \to B \not= 0$).
\end{proposition}

\subsection{Pairs of distinguished cubics}

\begin{figure}[htbp]
\centering 
\includegraphics{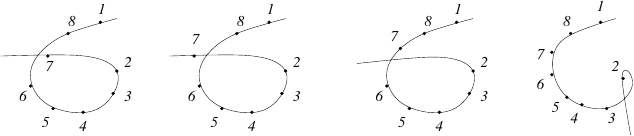}
\caption{\label{disting}   Cubics $(1-, 7), (7, 1-), (1-, 67)$ and $(1-, E)$}
\end{figure}

\begin{figure}[htbp]
\centering
\includegraphics{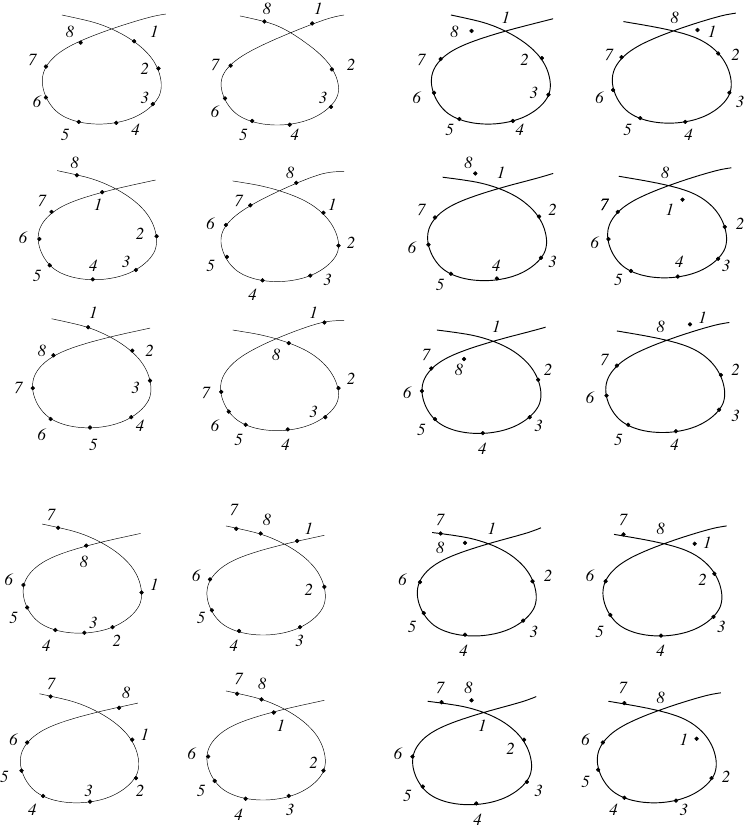}
\caption{\label{example}  Cubics obtained with the close pairs $(\hat 8 = 1\pm, \hat 1 = 8\pm)$, and $(\hat 8 = \hat 1 =7+)$}
\end{figure}

In what follows, we call cubic $C_3$ a combinatorial cubic $(C_3, 1, \dots 8)$ (the position of $9$ is not yet specified).
We define an encoding for distinguished cubics, that is consistant with the action of $D_8$. 
Let $(1-, N), N = 3, \dots 8$ be the cubic obtained from $(1-, N)_{nod}$, shifting $N$ away from the node $X$ onto the arc $X2$ of the loop, and let $(N, 1-), N = 3, \dots 7$ be the cubic obtained from $(1-, N)_{nod}$, shifting $N$ away from the node onto the arc $X1$ of the odd component, see Figure~\ref{disting}.
For $(P, Q) = (3, 4), (4, 5), \dots, (8, 1)$, denote by $(1-, PQ)$, the cubic obtained from $(1-, P)_{nod}$ (or $(1-, Q)_{nod}$) shifting the node $X$ away from $P$ (or $Q$) into the interior of the arc $PQ$.
Let $(1-, E)$ ($E$ stands for {\em end\/}) be the cubic that could be defined as $(1-, N)$, with $N = 2$ or as $(1-, PQ)$, with $(P, Q) = (2, 3)$. We chose the specific notation to avoid double notation for one single cubic type. 
The cubics encoded hereabove are called {\em cubics of the family $1-$\/}. 
At last, denote by $(81, L)$ and $(81, C)$ the first and the second cubic in Figure~\ref{example}. The letters $L$ and $C$ stand respectively for {\em loop\/} and {\em odd component\/}. 

In section~5 where we classify the pencils, we will consider combinatorial cubics $(C_3, 1, \dots 9)$. To encode them, we use the notation defined here for $(C_3, 1, \dots, 8)$, enhanced with the position of $9$.
The cubic $C_3 = (C_3, 1, \dots, 8)$ is divided into $9$ successive arcs by the points $1, \dots 8, X$. 
If $C_3 = (81, L)$, then $9$ lies on the arc $XX$.
If $C_3$ is of the family $N\epsilon$, percourse $C_3$ starting from $N$ in the direction $\epsilon$, and write which oriented arc contains $9$, see Figures~\ref{maxpen}-\ref{pentwo}.

Let $1, \dots 8$ realize a list $L$ and a pencil $\mathcal{P}$. Assume we know $L$ and not $\mathcal{P}$. We explain hereafter how to determine pairs of distinguished cubics of $\mathcal{P}$. The position of $9$ is not specified, it can sometimes be found using Bezout's theorem. 
Move $1, \dots 8$, shifting $8$ towards $1$ till $8$ reaches $1$ (if $8 \to 1 = 0$) or $8$ reaches the first of the six conics separating it from $1$ (if $8 \to 1 \not= 0$), see Proposition~4. The only other base point of the pencil that moves is $9$. {\em  Assume that in the motion, the eight points are never on a cubic with node at one of them\/}. This condition is always achieved if the list is non-nodal. 
More generally, if the condition is achieved for a given set of points $1, \dots 8$, if $8 \to 1 \not= 0$, we say that the corresponding elementary change is realizable starting from this set of points.
The mobile pencil $\mathcal{P}(1, \dots 8)$ is preserved all along, except in the end, when it becomes a singular pencil $\mathcal{P}_{sing}$. Let $C^0_3, C^1_3$ and $C^2_3$ be the three singular cubics involved, see section~1.2.  
A {\em close pair\/} is a pair $\hat A = B\epsilon$, $\hat B = A\epsilon$
or $\hat A = \hat B = N\epsilon$, with $A, B$ consecutive. Up to the action of $D_8$ there are seven close pairs:
$\hat 8 = 1+, \hat 1 = 8+$, $\hat 8 = \hat 1 = N+$, $N = 2, \dots 7$.

\begin{figure}[htbp]
\centering 
\includegraphics{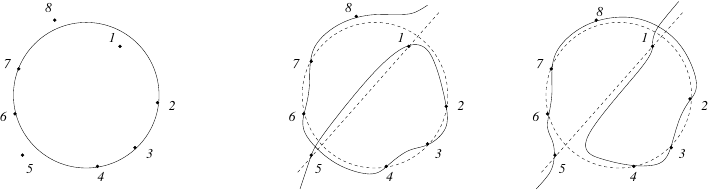}
\caption{\label{pair}  Cubics obtained with the elementary pair $(\hat 1: 5+, \hat 5: 8-)$}
\end{figure}

\begin{enumerate}
\item
If $8 \to 1 = 0$, the singular cubic $C^0_3$, with double point at $8 = 1$, is identical to both $(\hat 8, 1)$ and $(\hat 1, 8)$ (auxiliary cubics). 
The close pairs $\hat 8 = \hat 1 = 2+$ or $7-$ (deduced one from the other by $(+1)(48)$) are {\em inessential\/}: one does not know a priori whether the double points of the auxiliary cubics are isolated or not. Therefore, the cubics $C^1_3$, $C^2_3$ may be non-real. There are in total $16$ inessential close pairs (one orbit).
For the other close pair $(\hat 1, \hat 8)$ (essential close pairs), the auxiliary cubics have each a non-isolated double point. 
All along the motion $8 \to 1$, the position of $1$ with respect to $(\hat 1, 8)$ and the position of $8$ with respect to $(\hat 8, 1)$ are preserved. Using Bezout's theorem with these cubics, one finds out the corresponding pair $C^1_3$, $C^2_3$.
In all of the cases, any one of the four combinatorial data: cubic $(\hat 1, 8)$ enhanced with the position of $1$, cubic $(\hat 8, 1)$ enhanced with the position of $8$, cubic $C^1_3$, and cubic $C^2_3$ determines the other three.
The close pairs $(\hat 8 = 1+, \hat 1 = 8+)$ and $(\hat 8 = 1-, \hat 1 = 8-)$ give the same three admissible pairs $C^1_3$, $C^2_3$, they are shown in the upper part of Figure~\ref{example}, along with the auxiliary cubics.
\begin{center}
$(81, L), (81, C)$\\
$(8+, 1), (8-, 78)$\\
$(1+, 12), (1-, 8)$
\end{center}
The other close pairs $(\hat 8, \hat 1)$ give each rise to two admissible pairs $C^1_3$, $C^2_3$, see the tabular hereafter, where $N$ ranges from $3$ to $6$.
Note that all of these cubics $C^1_3$, $C^2_3$ are distinguished. The case $\hat 8 = \hat 1 = 7+$ is shown in the lower part of 
Figure~\ref{example}.
\begin{figure}[htbp]
\begin{tabular}{ c c c c }
$\hat 8 = \hat 1 = 7+$ & $\hat 8 = \hat 1 = 2-$ & $\hat 8 = \hat 1 = N+$ & $\hat 8 = \hat 1 = N-$\\
$(7+, 8), (1, 7+)$ & $(2-, 8), (12, C)$ & $(N+, 1), (8, N+)$ & $(N-, 8), (1, N-)$\\
$(78, C), (7+, 1)$ & $(2-, 1), (8, 2-)$ & $(N+, 8), (1, N+)$ & $(N-, 1), (8, N-)$\\
\end{tabular}
\end{figure}

\item
If $8 \to 1 \not= 0$, one has $C_3 = C_2 \cup (1P), P \in \{2, \dots, 7\}$. If: 
$P = 2$ and both $1$ and $2$ are outside of $C_2 = 34567$, then one does not 
know a priori whether the double points of the reducible cubic $C^0_3$ are 
real or complex conjugated. Therefore, the cubics $C^1_3$, $C^2_3$ may be 
non-real. Let $8$ goes to the outside of $C_2 = 34567$, $1, 2 > 34567$, the elementary change is: $\hat 1: 2- \to 2+$, $\hat 2: 1- \to 1+$, it is {\em inessential\/}. There are in total $16$ inessential elementary changes (one orbit).
In the other cases, either cubic $C^1_3$, $C^2_3$ may be obtained from the reducible cubic $C^0_3$ by perturbing one of the double points of $C^0_3$. All of the pairs obtained are distinguished cubics.
Indeed, let $s$ be the intersection of the line $(1P)$ with the interior of 
$C_2$, and let $s_1, s_2$ be the two arcs of $C_2$ on either side of $(1P)$. 
The loop of each cubic $C^1_3$, $C^2_3$ is obtained perturbing one of the 
$s_i \cup s, i \in \{1, 2\}$, and both $s_i \cup s$ contain some points 
among $1, \dots, 8$. 
We give an example in  Figure~\ref{pair}, the elementary pair is $(\hat 1: 5+, \hat 5: 8-)$, the corresponding pair of distinguished cubics is 
$(5, 8-), (5+, 1)$.
\end{enumerate}

The first elementary change $\hat 1: 2- \to 2+$, $\hat 2: 1- \to 1+$ in Figures~\ref{elem}-\ref{repres} is inessential. 
The changes of distinguished cubics corresponding to the last $18$ elementary changes are displayed in Figure~\ref{changes}. We say that these changes are {\em essential\/}. An elementary pair appearing in some (in)essential elementary change will be called (in)essential elementary pair. Up to the action of $D_8$, two pairs are both elementary and close: $(\hat 1 = 8+, \hat 8 = 1+)$ (essential only as close pair), and $(\hat 1 = \hat 8 = 2+)$ (essential only as elementary pair).  
If a list $L$ is non-nodal, each essential elementary pair of $L$
gives a pair of distinguished cubics, see Figure~\ref{changes}, and each essential close pair of $L$ gives two or three admissible pairs of distinguished cubics. The rest of the list allows to determine which pair of cubics is actually realized.

\begin{figure}
\begin{tabular}{c c c c c}
$(1+, 2) \to (1-, E)$ & & $(1-, E) \to (1+, 2)$ & & $(8-, 2) \to (23, L)$\\
$(1, 3+) \to (1+, 12)$ & & $(1+, 12) \to (1, 3+)$ & & $(81, L) \to (3+, 1)$\\
 & & \\
$(3, 1+) \to (3, 1-)$ & & $(3, 1-) \to (3, 1+)$ & & $(1+, 3) \to (1-, 3)$\\
$(1, 3-) \to (1, 3+)$ & & $(1, 3+) \to (1, 3-)$ & & $(12, C) \to (1, 4+)$\\
 & & \\
$(1-, 3) \to (1+, 3)$ & & $(8-, 3) \to (2+, 3)$ & & $(2+, 3) \to (8-, 3)$\\
$(1, 4+) \to (12, C)$ & & $(4+, 1) \to (2-, 1)$ & & $(2-, 1) \to (4+, 1)$\\
 & & \\
$(4, 1+) \to (4, 1-)$ & & $(4, 1-) \to (4, 1+)$ & & $(1+, 4) \to (1-, 4)$\\
$(1, 4+) \to (1, 4-)$ & & $(1, 4-) \to (1, 4+)$ & & $(1, 5+) \to (1, 3-)$\\
 & & \\
$(1-, 4) \to (1+, 4)$ & & $(8-, 4) \to (2+, 4)$ & & $(2+, 4) \to (8-, 4)$\\
$(1, 3-) \to (1, 5+)$ & & $(3-, 1) \to (5+, 1)$ & & $(5+, 1) \to (3-, 1)$\\
 & & \\
$(5, 1+) \to (5, 1-)$ & & $(1+, 5) \to (1-, 5)$ & & $(8-, 5) \to (2+, 5)$\\
$(1, 5-) \to (1, 5+)$ & & $(1, 4-) \to (1, 6+)$ & & $(6+, 1) \to (4-, 1)$\\
\end{tabular}
\caption{\label{changes} Changes of pairs of distinguished cubics}
\end{figure}

\subsection{Simultaneous changes of lists and pencils, exceptions}

Let $1, \dots, 8$ be eight generic points lying in strictly convex position in $\mathbb{R}P^2$, realizing a list $L(1, \dots 8)$ and a pencil $\mathcal{P}=\mathcal{P}(1, \dots, 8)$. 
Move the eight points, keeping them distinct and strictly convex.
At some moment, a degeneration $\mathcal{P} \to \mathcal{P}_{sing}$ occurs, after which one gets a new pencil $\mathcal{P}'$, see section 1.2. We are interested only in the degeneration that affect the combinatorial pencil, so we leave aside the cases 1b, 2b, 2c, 2d.
When the list undergoes an elementary change, six of the eight chosen base points become coconic, $9$ becomes aligned with the other two. If the change is essential, the new pencil $\mathcal{P}'$ is obtained from $\mathcal{P}$ by replacing the pair of distinguished cubics $C^1_3$, $C^2_3$ by a new pair $C^3_3$, $C^4_3$ (case 2a).
If the change is inessential, we are in case 2a or 2b. 
In case 2a, one only of the two cubics $C^1_3$, $C^2_3$ is distinguished, and one recovers the same pair of combinatorial cubics after the degeneration (left to the reader). 
Performing an inessential change on the list doesn't change the pencil.
Conversely, one may sometimes change the pencil without changing the list. The degeneration occurs letting one of the first eight base points $A$ come together with $9$ (1a). 
The cubic $C_3^0$ of $\mathcal{P}_{sing}$ passes through $1, \dots 8$ and has a node at one of these points, the corresponding list $L(1, \dots 8)$ is nodal. 
Move further, the new generic pencil $\mathcal{P}'$ is deduced from 
$\mathcal{P}$ swapping the positions of $A$ and $9$ on all of the eight distinguished (combinatorial) cubics. 
Performing an inessential change on a nodal list yields again a nodal list.
So, if a pencil is realizable by two different lists, these lists are both non-nodal or both nodal. In the latter case, we say that the pencil is {\em nodal\/}. 


\section{Lists with configurations of six coconic points, pencils
with reducible cubics}

\subsection{Lists obtained perturbing four reducible cubics}

Among the principal lists, exactly ten have a set of four disjoint elementary pairs.
The lists $L_n$ with $n = 11, 12, 13, 14, 19, 20, 21, 22$, have each
disjoint elementary pairs $(\hat 1, \hat 8)$, $(\hat 2, \hat 7)$, $(\hat 3, \hat 6)$, $(\hat 4, \hat 5)$.
Each of the two lists $L_n$ with $n = 88, 92$ has disjoint elementary pairs 
$(\hat 1, \hat 5)$, $(\hat 2, \hat 6)$, $(\hat 3, \hat 7)$, $(\hat 4, \hat 8)$. 
We explain hereafter how to realize these lists directly.

\begin{figure}[htbp]
\centering 
\includegraphics{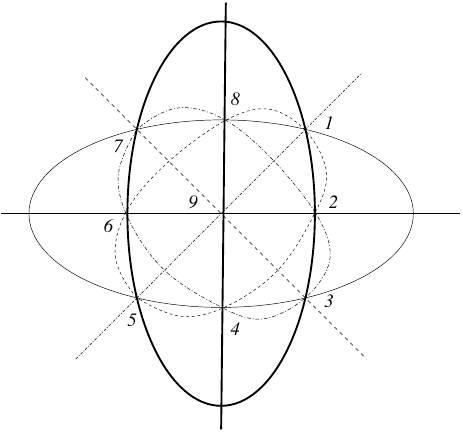}
\caption{\label{list1} First configuration of points with four reducible cubics}
\end{figure}

Consider a pair of ellipses intersecting at four points, $1, 3, 5, 7$, see  Figure~\ref{list1}. 
Denote by $9$ the intersection of the diagonal lines $15$ and $37$. Draw a vertical and a horizontal line, both passing through $9$.
Let $4, 8, 2, 6$ be four supplementary points chosen such as: $4, 8$ are the intersections of the vertical line with one ellipse; $2, 6$ are the intersections of the horizontal line with the other ellipse; each pair of points lying on one ellipse is in the interior of the other. By construction, the points $1, \dots 8$ lie in convex position. Moreover, there exist two supplementary conics passing through six points: $234678$ and $124568$.
As a matter of fact, the pencil of cubics $\mathcal{P}$ determined by $1, \dots, 8$ has $9$ as ninth base point and four distinguished cubics, all of them reducible: $498 \cup 123567$, $296 \cup 134578$, $397 \cup 124568$ and $195 \cup 234678$.
Let us perturb the pencil, moving slightly the points $8$ and $9$. As $1, 2, 3, 5, 6, 7$ are on a conic, $8, 4$ and $9$ must stay aligned.
The point $8$ leaves the three conics $23467, 12456, 13457$. 
Letting $8$ cross each conic from the inside to the outside yields the following elementary moves:
Conic $23467$: $(\hat 1 = 5-, \hat 5 = 1+) \to (\hat 1 = 5+, \hat 5 = 1-)$. 
Conic $12456$: $(\hat 3 = 7-, \hat 7 = 3+) \to (\hat 3 = 7+, \hat 7 = 3-)$. 
Conic $13457$: $(\hat 2 = 7+, \hat 6 = 1-) \to (\hat 2 = 5-, \hat 6 = 3+)$.
Let us now move $3$ away from the conic $12567$. 
Letting $3$ cross $12567$ from the inside to the outside yields the following elementary move:
$(\hat 4 = 7-, \hat 8 = 5+) \to (\hat 4 = 1+, \hat 8 = 3-)$.
The first perturbation may be done so as to realize six different positions of $8$ with respect to the set of conics
$23467$, $12456$, $13457$. Then, move $3$ to the outside of $12567$. We obtain the first five lists, and the last list of
Figure~\ref{lsix}, among which we find $L_{88}$ and $L_{92}$.

\begin{figure}[htbp]
\centering 
\includegraphics{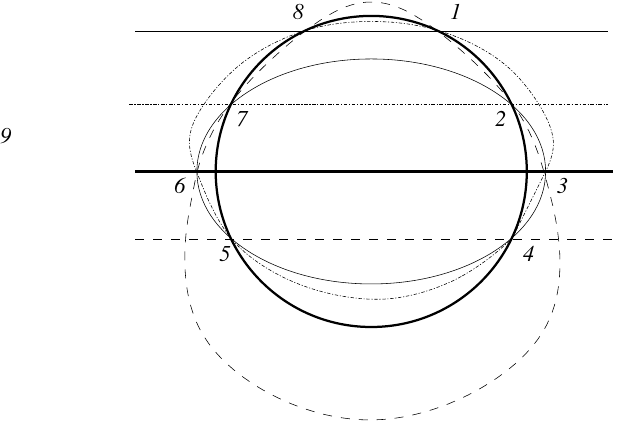}
\caption{\label{list2} Second configuration of points with four reducible cubics}
\end{figure}

Consider a pair of ellipses intersecting at four points $2, 4, 5, 7$, see Figure~\ref{list2}. 
Denote by $9$ the intersection of the lines $27$ and $45$. Draw a line $\Delta$ passing through $9$ and cutting the ellipses on their arcs $57$ and $24$. Let $6$ and $3$ be the intersections of $\Delta$ with one ellipse, chosen so that these points lie outside of the second ellipse. Draw a line $\Delta'$ passing through $9$ and cutting the second ellipse at two points $8, 1$ on the
arc $72$. One may choose $\Delta'$ in such a way that the points $1, \dots 8$ lie in convex position.
By construction, there exist two supplementary conics passing through six of the points: $123678$ and $134568$.
The pencil of cubics $\mathcal{P}$ determined by $1, \dots, 8$ has $9$ as ninth base point and four distinguished cubics, all of them reducible: $189 \cup 234567$, $369 \cup 124578$, $459 \cup 123678$ and $279 \cup 134568$.
Let us perturb the pencil, moving the points $8$ and $9$.
The point $8$ leaves the three conics $13456$, $12457$, $12367$.
Letting $8$ cross each conic from the inside to the outside yields the following elementary moves:
Conic $13456$: $(\hat 2 = 6-, \hat 7 = 1-) \to (\hat 2 = 8+, \hat 7 = 3+)$.
Conic $12457$: $(\hat 3 = 6+, \hat 6 = 3+) \to (\hat 3 = 6-, \hat 6 = 3-)$.
Conic $12367$: $(\hat 4 = 3-, \hat 5 = 3-) \to (\hat 4 = 6+, \hat 5 = 6+)$.
Let us now move $7$ away from the conic $23456$. 
Letting $7$ cross $23456$ from the inside to the outside yields the following elementary move:
$(\hat 1 = 8+, \hat 8 = 1+) \to (\hat 1 = 8-, \hat 8 = 1-)$.
The first perturbation may be done so as to realize six different positions of $8$ with respect to the set of conics
$13456$, $12457$, $12367$. Then, move $7$ to the inside of $23456$. We obtain the lists $L_{11}, L_{12}, L_{14}, L_{19}, L_{21}, L_{22}$.
Starting again from the pencil with four reducible cubics we realize the two missing lists as follows.
List $L_{13}$: move first $7$ to the outside of $12458$ and the inside of $12368$; then, move $8$ to the outside of
$13456$. List $L_{20}$: move first $6$ to the outside of $23457$ and the inside of $12378$; then move $8$ to the outside of
$12457$.

\begin{proposition}
Let $1, \dots 8$ be eight points lying in convex position in the plane and let $k$ be the number of conics passing through exactly six of them. One has $k \leq 4$. If $k = 4$, then
the points realize, up to the action of $D_8$, one of the two non-generic lists shown in Figures~ \ref{list1}-\ref{list2}.  
The orbit of the first list has two elements, the orbit of the second list has eight elements.
\end{proposition}

\begin{description}
\item[Proof:]
Perturbing slightly the configuration $1, \dots 8$ must yield a generic list with four distinct elementary pairs, otherwise stated a list that is in the orbit of some of the $10$ lists $L_n$. 
Up to the action of $D_8$, the original configuration $1, \dots 8$ must be as shown in Figure~\ref{list1} or \ref{list2}. The list $l$ of Figure~\ref{list1} is encoded by the data: $2, 6 < 134578$, $4, 8 < 123567$, $1, 5 > 234678$, $3, 7 > 124568$.
This list  is invariant by the action of $id$, $15$, $26$, $37$, $48$, $+2$, $-2$, $+4$.
The list $l'$ of Figure~\ref{list2} is encoded by the data: $1, 8 > 234567$, $5, 4 < 123678$, $2, 7 < 134568$, $3, 6 > 124578$. This list is invariant by the action of $(+1)(48)$.
$\Box$
\end{description}

A pencil of cubics with eight base points lying in convex position in the real plane (no seven of them being coconic) has at most four reducible cubics, the corresponding four lines pass all through the ninth base point. In the next section, we drop the condition of convexity and search for more singular pencils.

\subsection{A singular pencil with base points in non-convex position}

Let us say that a cubic is  {\em completely reducible\/} if it is the product of three lines.
A complex pencil contains at most four such cubics \cite{st}, \cite{yu}, the corresponding twelve lines and the nine base points are such that: each point lies on four lines and each line passes through three points.  Recall that the nine inflection points of a complex cubic $C_3$ realize such a configuration. The {\em Hessian pencil\/} associated to $C_3$ is the pencil 
generated by $C_3$ and its Hessian, based at the inflection points of $C_3$. 
The Hessian pencils realize the upper bound of four completely reducible cubics.

\begin{figure}[htbp]
\centering 
\includegraphics{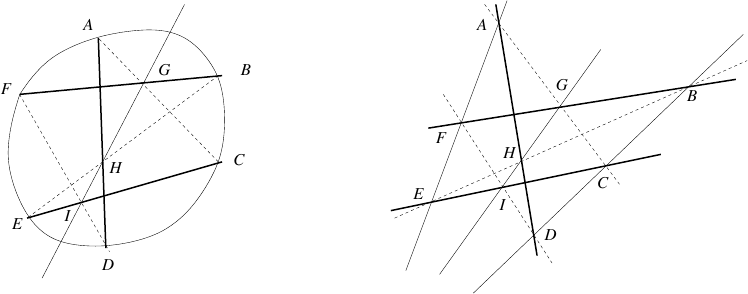}
\caption{\label{pascal}  Pascal's and Pappus' theorems}
\end{figure}

Let us now go back to pencils with only real base points.
A pencil with nine real base points cannot have four completely reducible cubics:  
the Sylvester-Gallai theorem states that given $n$ points in the real plane, they are either all collinear or there exists a line containing exactly two of them (see e.g.\cite{B-M}) . 
Consider the pencil generated by the two completely reducible cubics (one bold, one dotted) in Figure~\ref{pascal}, the base points are denoted by $A, \dots I$. We recover an elementary proof of Pascal's theorem: assume that $G, H, I$ are on a line, and let $C_2$ be a conic passing through five of the other base points. Then, the cubic $(GHI) \cup C_2$ belongs to the pencil, otherwise stated, $A, B, C, D, E, F$ are coconic. Conversely, if $A, \dots F$ are coconic, then $G, H, I$ are aligned. The particular case where $C_2$ is the product of two lines is Pappus' theorem, the pencil here has three completely reducible cubics. 

\begin{figure}[htbp]
\centering 
\includegraphics{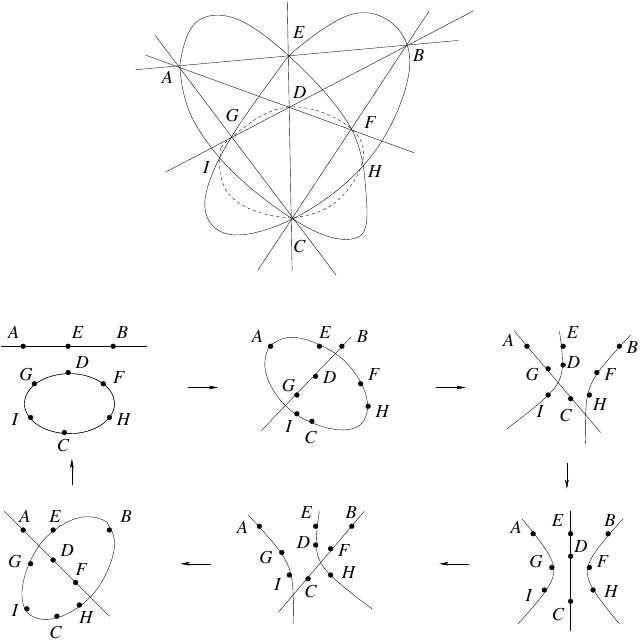}
\caption{\label{sixcub}  Pencil with six reducible cubics}
\end{figure}

Let us now search for a pencil having the maximal number of six reducible cubics. Each of the corresponding six lines must pass through exactly three base points. If one point lies on four of these lines, then a fifth line intersects these four at four base points, contradiction. If each base point lies on at most two lines, draw five lines, each of them passes through four base points, contradiction again.  Thus, one point lies on exactly three lines, six other base points are distributed pairwise on these lines. Draw two further lines, they must intersect at one of these
six points, otherwise they would pass each through four base points. We get thus the following distribution of the base points on the lines and the conics: four points $A, B, C, D$ lie each on three lines and three conics, three points $D, E, F$ lie each on two lines and four conics and the other two points $H, I$ lie each on six conics. Assume that $H, I$ are also real. Up to the action of the symmetric group $\mathcal{S}_4$ on the first chosen base points $A, B, C, D$, the sequence of reducible cubics is as shown in Figure~\ref{sixcub}. Note that the non-singular cubics of the pencil are all disconnected. 

\subsection{Symmetric lists}

We proved in section 2.2 that a generic list is preserved by no element of $D_8$, otherwise stated, each orbit contains $16$ lists. Let us now search for the non-generic lists that are invariant for some element of $D_8$, we met two of them, $l$ and $l'$, in section 4.2. Let $\lambda_0$ be the list $12345678$ consisting of eight coconic points. We say that a list is {\em almost generic\/} if it has one unique configuration of six coconic points.

\begin{proposition}
Each element of $D_8$ leaves some non-generic lists invariant:
$(\pm 1)$ and $(\pm 3)$: $\lambda_0$;
$(\pm 2)$ and $(\pm 4)$: $\lambda_0$, $l$ and $(+1)l$;
$(+1)(48)$: $\lambda_0$, $l'$, $(+4)l' = (+1)(26)l'$;
$15$: $\lambda_0$, four lists with seven coconic points, two lists with each four sets of six coconic points, ten lists with each three sets of six coconic points, two lists with each two sets of six coconic points and $24$ almost generic lists, divided in two groups of $12$ that are mapped one onto the other by $(+4)$ (or equivalently $37$).  
All of these lists have a conic $234678$. See description ahead and Figures~\ref{alone}-\ref{althree}.   
\end{proposition}

In what follows, $a$ stands for a cyclic permutation, and $\sigma$ for a symmetry $15$, $26$, $37$ or $48$. Recall that any two cyclic permutations commute, whereas 
$a\sigma = \sigma (-a)$. Let $m$ be a cyclic permutation, and $L$ be a list such that $m\cdot L = L$. Any list $L'$ in the orbit of $L$ is also invariant under the action of $m$ and $-m$. Up to the action of $D_8$, we may assume that $L$ contains a conic $C_2 = 123456$, $123457$, $123467$ or $123567$. Moreover, $L$ contains the images of this conic by $m$, $m \circ m$ \dots and if $P < C_2$, then $m \cdot P < m \cdot C_2$.  One finds out easily that $(\pm 1)$, $(\pm3)$ preserve only $\lambda_0$, whereas $(\pm 2)$, $(\pm 4)$ preserve each three lists: $\lambda_0$, $l$ and $(+1)l$.
A list invariant by $(+1)(48)$ must contain the conics $124578$, $134568$, $123678$ and 
$234567$. Thus, $(+1)(48)$ preserves three lists: $\lambda_0$, $l'$ and $(+4)l' = (+1)(26)l'$.
A list invariant by $15$ must contain a conic $234678$. Let us first look for the almost generic lists. 
Let $L$ be such a list, the orbit of $L$ contains eight elements, each symmetry $\sigma$ leaves two of them invariant. If $\sigma \cdot L' = L'$, one may perturb $L'$ in two ways to get generic lists that are deduced from one another by $\sigma$. The $16$ generic lists thus obtained form an orbit of $D_8$. We select all of the principal lists having an elementary pair $(\hat N, \hat M)$, with $(N, M)$ = $(1, 5)$, $(2, 6)$, $(3, 7)$ or $(4, 8)$, and such that the non-generic list obtained making the six points different from $N, M$ coconic is invariant by the symmetry $NM$.  Up to cyclic permutations, we get thus $24$ almost generic lists invariant by $15$, they split in two groups of $12$ that are deduced from one another by the action of $+4$ (or $37$). These lists may be obtained also as perturbations of more singular lists invariant by $15$.


\begin{enumerate}
\item
If $1$ or $5$ lies on $234678$
\begin{enumerate}
\item
One list with eight coconic points: $\lambda_0 = 12345678$ 
\item
Four lists with seven coconic points: $5 < 1234678$, $5> 1234678$, 
$1 < 2345678$, $1 > 2345678$.
\end{enumerate}
\item
If $1 > 234678$ and $5 < 234678$
\begin{enumerate}
\item
Four almost generic lists $l_1, l_2, l_3, l_4$ 
\item
Three lists $\lambda_1, \lambda_2, \lambda_3$ having each two more configurations of six coconic points:  
\begin{center}
$\lambda_1: 123456$, $145678$; \qquad $\lambda_2: 123457$, $135678$; \\ 
$\lambda_3: 123458$, $125678$
\end{center}
\end{enumerate}
\item
If $1 < 234678$ and $5 > 234678$
\begin{enumerate}
\item
Four almost generic lists $l'_i = (+4)(l_i)$
\item
Three lists $\lambda'_i = (+4)(\lambda_i)$ having each two more configurations of six coconic points: 
\begin{center}
$\lambda'_1: 123458$, $125678$; \qquad $\lambda'_2: 123457$, $135678$; \\
$\lambda'_3: 123456$, $145678$
\end{center}
\end{enumerate}
\item
If $1 > 234678$ and $5 > 234678$
\begin{enumerate}
\item
Eight almost generic lists $l_5, l_6, l_7, l_8$, $l'_5, l'_6, l'_7, l'_8$, 
\item
Four lists having each one more configuration of six coconic points: 
\begin{center}
$\lambda_4: 123578$; \qquad $\lambda_5 = (+4)(L_4): 134567$; \qquad $\lambda_6, \lambda_7 = (+4)(L_6): 124568$ 
\end{center}
\item
Two lists having each two more configurations of six coconic points:
\begin{center}
$\lambda_8, \lambda_9 = (+4)(\lambda_8)$: $123567, 134578$
\end{center}
\item
One list with three more configurations of six coconic points:
\begin{center}
$\lambda_{10} = l: 124568, 134578, 123567$
\end{center}
\end{enumerate}
\item
$1 < 234678$ and $5 < 234678$
\begin{enumerate}
\item
Eight almost generic lists $l_9, l_{10}, l_{11}, l_{12}$, $l'_9, l'_{10}, l'_{11}, l'_{12}$
\item
Four lists having each one more configuration of six coconic points: 
\begin{center}
$\lambda_{11}: 123578$; \qquad $\lambda_{12} = (+4)(\lambda_{11}): 134567$; \qquad $\lambda_{13}, \lambda_{14} = (+4)(\lambda_{13}): 124568$ 
\end{center}
\item
Two lists having each two more configurations of six coconic points:
\begin{center}
$\lambda_{15}, \lambda_{16} = (+4)(\lambda_{15}): 123567, 134578$
\end{center}
\item
One list having three more configurations of six coconic points:
\begin{center}
$\lambda_{17} = (+1)(\lambda_{10}) = (+1)(l): 124568, 123567, 134578$
\end{center}
\end{enumerate}
\end{enumerate}

\begin{figure}
\begin{tabular}{ c c c c c }
 & $l_1$ & $l_2$ & $l_3$ & $l_4$\\
$\hat 1$ & $6+ \leftrightarrow 4-$ & $6+ \leftrightarrow 4-$ & $6+ \leftrightarrow 4-$ & $6+ \leftrightarrow 4-$\\
$\hat 2$ & $1-$ & $4-$ & $4+$ & $4-$\\
$\hat 3$ & $1-$ & $4+$ & $4+$ & $4-$\\
$\hat 4$ & $1-$ & $3+$ & $1-$ & $3-$\\
$\hat 5$ & $1- \leftrightarrow 1+$ & $1- \leftrightarrow 1+$ & $1- \leftrightarrow 1+$ & $1- \leftrightarrow 1+$\\
$\hat 6$ & $1+$ & $7-$ & $1+$ & $7+$\\
$\hat 7$ & $1+$ & $6-$ & $6-$ & $6+$\\
$\hat 8$ & $1+$ & $6+$ & $6-$ & $6+$\\
 & & & & \\
 & $l'_1$ & $l'_2$ & $l'_3$ & $l'_4$\\
$\hat 1$ & $5- \leftrightarrow 5+$ & $5- \leftrightarrow 5+$ & $5- \leftrightarrow 5+$ & $5- \leftrightarrow 5+$\\
$\hat 2$ & $5+$ & $3-$ & $5+$ & $3+$\\
$\hat 3$ & $5+$ & $2-$ & $2-$ & $2+$\\
$\hat 4$ & $5+$ & $2+$ & $2-$ & $2+$\\
$\hat 5$ & $2+ \leftrightarrow 8-$ & $2+ \leftrightarrow 8-$ & $2+ \leftrightarrow 8-$ & $2+ \leftrightarrow 8-$\\
$\hat 6$ & $5-$ & $8-$ & $8+$ & $8-$\\
$\hat 7$ & $5-$ & $8+$ & $8+$ & $8-$\\
$\hat 8$ & $5-$ & $7+$ & $5-$ & $7-$\\
\end{tabular}
\caption{\label{alone} Almost generic symmetric lists (first)}
\end{figure}

\begin{figure}
\begin{tabular}{ c c c c c }
 & $l_5$ & $l_6$ & $l_7$ & $l_8$\\
$\hat 1$ & $5- \leftrightarrow 5+$ & $5- \leftrightarrow 5+$ & $5- \leftrightarrow 5+$ & $5- \leftrightarrow 5+$\\
$\hat 2$ & $1+$ & $7+$ & $7+$ & $7-$\\
$\hat 3$ & $1+$ & $7-$ & $7+$ & $1+$\\
$\hat 4$ & $1+$ & $1+$ & $1+$ & $1+$\\
$\hat 5$ & $1+ \leftrightarrow 1-$ & $1+ \leftrightarrow 1-$ & $1+ \leftrightarrow 1-$ & $1+ \leftrightarrow 1-$\\
$\hat 6$ & $1-$ & $1-$ & $1-$ & $1-$\\
$\hat 7$ & $1-$ & $3+$ & $3-$ & $1-$\\
$\hat 8$ & $1-$ & $3-$ & $3-$ & $3+$\\
 & & & & \\
 & $l'_5$ & $l'_6$ & $l'_7$ & $l'_8$\\
$\hat 1$ & $5- \leftrightarrow 5+$ & $5- \leftrightarrow 5+$ & $5- \leftrightarrow 5+$ & $5- \leftrightarrow 5+$\\
$\hat 2$ & $5-$ & $5-$ & $5-$ & $5-$\\
$\hat 3$ & $5-$ & $7+$ & $7-$ & $5-$\\
$\hat 4$ & $5-$ & $7-$ & $7-$ & $7+$\\
$\hat 5$ & $1+ \leftrightarrow 1-$ & $1+ \leftrightarrow 1-$ & $1+ \leftrightarrow 1-$ & $1+ \leftrightarrow 1-$\\
$\hat 6$ & $5+$ & $3+$ & $3+$ & $3-$\\
$\hat 7$ & $5+$ & $3-$ & $3+$ & $5+$\\
$\hat 8$ & $5+$ & $5+$ & $5+$ & $5+$\\
\end{tabular}
\caption{\label{altwo} Almost generic symmetric lists (second)}
\end{figure}

\begin{figure}
\begin{tabular}{ c c c c c }
 & $l_9$ & $l_{10}$ & $l_{11}$ & $l_{12}$\\
$\hat 1$ & $6+ \leftrightarrow 4-$ & $6+ \leftrightarrow 4-$ & $6+ \leftrightarrow 4-$ & $6+ \leftrightarrow 4-$\\
$\hat 2$ & $8+$ & $6-$ & $6+$ & $8-$\\
$\hat 3$ & $8-$ & $8+$ & $8+$ & $8-$\\
$\hat 4$ & $8-$ & $8-$ & $8+$ & $8-$\\
$\hat 5$ & $8- \leftrightarrow 2+$ & $8- \leftrightarrow 2+$ & $8- \leftrightarrow 2+$ & $8- \leftrightarrow 2+$\\
$\hat 6$ & $2+$ & $2+$ & $2-$ & $2+$\\
$\hat 7$ & $2+$ & $2-$ & $2-$ & $2+$\\
$\hat 8$ & $2-$ & $4+$ & $4-$ & $2+$\\
 & & & & \\
 & $l'_9$ & $l'_{10}$ & $l'_{11}$ & $l'_{12}$\\
$\hat 1$ & $6+ \leftrightarrow 4-$ & $6+ \leftrightarrow 4-$ & $6+ \leftrightarrow 4-$ & $6+ \leftrightarrow 4-$\\
$\hat 2$ & $6+$ & $6+$ & $6-$ & $6+$\\
$\hat 3$ & $6+$ & $6-$ & $6-$ & $6+$\\
$\hat 4$ & $6-$ & $8+$ & $8-$ & $6+$\\
$\hat 5$ & $8- \leftrightarrow 2+$ & $8- \leftrightarrow 2+$ & $8- \leftrightarrow 2+$ & $8- \leftrightarrow 2+$\\
$\hat 6$ & $4+$ & $2-$ & $2+$ & $4-$\\
$\hat 7$ & $4-$ & $4+$ & $4+$ & $4-$\\
$\hat 8$ & $4-$ & $4-$ & $4+$ & $4-$\\
\end{tabular}
\caption{\label{althree} Almost generic symmetric lists (third)}
\end{figure}

\newpage
\section{Classification of the pencils of cubics}

\subsection{Nodal pencils}

\begin{figure}[htbp]
\centering 
\includegraphics{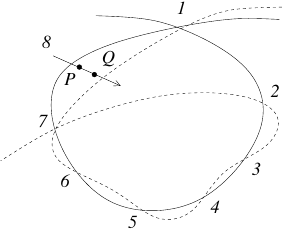}
\caption{\label{path} The path $p$}
\end{figure}
\begin{figure}[htbp]
\centering 
\includegraphics{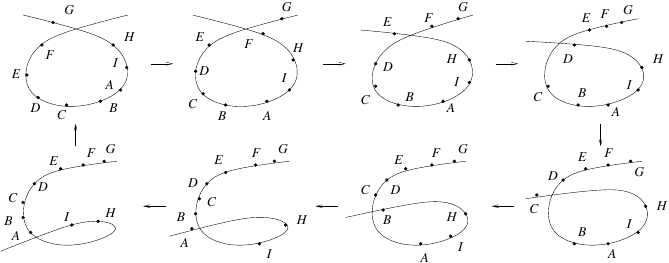}
\caption{\label{pen} Nodal pencils obtained from the list $\max(\hat 1 = 8-)$}
\end{figure}

We will use the method exposed in section~3.2 to construct the nodal pencils.
Let $T$ be the triangle $(12), (67), (17)$ containing $8$.
The condition that $1, \dots, 8$ realizes the list $\max (\hat 1 = 8-)$ splits into eight disjoint subconditions. There is an ordering $14567 > (\hat 8, 1) > (\hat 8, 7) > (\hat 8, 6) > \dots > (\hat 8, 2) > 34567$ in $T$:
the point $8$ lies between two consecutive of these curves. When $8$ lies on a cubic $(\hat 8, N)$, otherwise stated, when the eight chosen base points are on a cubic $(1-, N)_{nod}$, the pencil is singular, with $9 = N$.
By Bezout's theorem with the cubics $(\hat 8, N)$, the degeneration $9 = 8$ may occur only if $8$ lies between $(\hat 8, 1)$ and $(\hat 8, 7)$.
Start with $8$ between $14567$ and $(\hat 8, 1)$, the elementary change 
$8$ goes to the outside of $14567$, $2, 3 < 14567$ writes:
$\hat 2: 1- \to 4+, \hat 3: 1- \to 4+$, it is realizable. The essential elementary pair $(\hat 2 = 1-, \hat 3 = 1-)$ gives thus two cubics: $(1-, 3)$, $(12, L)$. The ninth base point $9$ must lie on the odd component of $(12, L)$, by Bezout's theorem, it lies then on the arc $1X$ of $(1-, 3)$. We have $3 \to 4 = 0$, the close pair $(\hat 3 = 1-, \hat 4 = 1-)$ gives two admissible pairs of cubics: $(3, 1-), (1-, 4)$ or $(1-, 3), (4, 1-)$. We know already that $(1-, 3)$ belongs to the pencil, so the correct pair is the first one. Using the other pairs of consecutive points $A, B$ with $A \to B = 0$, we get other cubics: the pair $4, 5$ gives again $(4, 1-)$ and a new cubic $(1-, 5)$, the pairs $5, 6$ and $6, 7$ give $(6, 1-)$ and $(1-, 7)$. Note that $2, 3$ gives $(\hat 2 = 1-, \hat 3 = 1-)$, inessential as close pair.
The position of $9$ on each cubic is again obtained with Bezout's theorem. 
Up to now, we have six cubics. Note that all of them have an arc $81$, so the missing cubics correspond to the openings of this arc. Starting from $(12, L)$, move into the portion formed of cubics with ovals. The next distinguished cubic will be
$(81, L)$, the overnext is $(81, C)$. The complete pencil is shown in the first row of
Figure~\ref{maxpen}. 
Let now $8$ percourse a path $p$ crossing successively $(\hat 8, 1)$ and $(\hat 8, 7)$. We will prove that at some moment, $8$ must indeed come together with $9$. Letting $8$ cross $(\hat 8, 1)$ changes the pencil swapping the positions of
$1$ and $9$ on all combinatorial cubics, see second row of Figure~\ref{maxpen}. 
The point $8$ lies now in the zone between the arcs $17$ of $(\hat 8, 1)$ and $(\hat 8, 7)$, this zone is divided in an upper and a lower part by the path $p$. Applying Bezout's theorem with $(\hat 8, 1)$ and $(1-, 81), 2X$ on one hand, and
with $(\hat 8, 7)$ and $(1-, 7), 18$ on the other hand, we get that $9$ lies also in this zone. When $8$ percourses the path $p$ from $(\hat 8, 1)$ and $(\hat 8, 7)$, the point $9$ moves from $1$ to $7$. Let $9$ cross $p$ at some point $P$ and assume that at this moment, $8$ has the position $Q$. Assume first that $Q$ lies ahead from $P$, see Figure~\ref{path}. When $8$ had previously the position $P$, the point $9$ was in still in the upper zone. But on the other hand, $9$ must have been positioned at $Q$. This is a contradiction. One gets a similar contradiction taking $Q$ behind $P$ on the path. So, at some moment, one must have $9 = 8$.
Letting $8$ cross successively $(\hat 8, 1), \dots  (\hat 8, 2)$, we obtain in total the first nine pencils of Figure~\ref{maxpen}, they are drawn in Figure~\ref{pen} where

$(G, H, I, A, B, C, D, E, F) = (9, 2, 3, 4, 5, 6, 7, 8, 1)$, $(1, 2, 3, 4, 5, 6, 7, 8, 9)$, 

$(1, 2, 3, 4, 5, 6, 7, 9, 8)$, $(1, 2, 3, 4, 5, 6, 9, 7, 8)$, $(1, 2, 3, 4, 5, 9, 6, 7, 8)$,

$(1, 2, 3, 4, 9, 5, 6, 7, 8)$, $(1, 2, 3, 9, 4, 5, 6, 7, 8)$, $(1, 2, 9, 3, 4, 5, 6, 7, 8)$, 

$(1, 9, 2, 3, 4, 5, 6, 7, 8)$.

If $8$ lies outside of the loop of $(\hat 8, 1)$, then one gets the pencil with $G = 9$.
If $8$ lies between $(\hat 8, 1)$ and $(\hat 8, 7)$, one gets the next two pencils with $(E, F) = (8, 9)$ and $(9, 8)$.
The other positions of $8$ give rise to the other pencils, switching successively $9$ with $7, 6, \dots, 2$.




Let $1, \dots 8$ realize $\max(\hat 1 = 8-)$.
The inessential elementary change $\hat 1: 8- \to 8+$, $\hat 8: 1- \to 1+$ ($7$ enters the conic $23456$, $1, 8$ outside of $23456$) is realizable if and only if the points realize one of the first three nodal pencils.
Let $1, \dots 8$ realize one of these pencils, we may move $8$ leaving the other points fixed, till the eight points lie on a cubic $(1\pm, 1)_{nod}$ or $(1-, 8)_{nod}$, see Figure~\ref{combi}. The conic $23456$ has a sixth intersection point $P$ with the loop of the cubic.
The point $P$ lies on some arc $AB$ where $A$ and $B$ are two consecutive points among $1, \dots 7$. If we move $B$ towards $A$ on the cubic, $P$ moves in the opposite direction, at some moment one has $P = B$ (the cubic and the conic have an ordinary tangency point), then the positions of the two points are swapped.   
So, up to some swaps like this, we may arrange that $P$ lies on the arc $67$ of the cubic. 
Move $7$ towards $6$ along the cubic till $7$ crosses $23456$, we are done. Let now $1, \dots 8$ realize one of the other six pencils. The  sixth intersection point of $23456$ with the cubic $(7, 1-)$ is on the loop of this cubic, whereas $7$ is on the odd component. The elementary change is not realizable. 

The elementary change $\hat 1: 8- \to 8+$, $\hat 8: 1- \to 1+$ replaces $\max(\hat 1 = 8-)$ by $\max(\hat 1 = 8+)$, and leaves the combinatorial pencil unchanged.
The nodal list $\max(\hat 1 = 8+)$ is thus realizable by the three pencils of Figure~\ref{pen}  with 

$(G, H, I, A, B, C, D, E, F) = (9, 2, 3, 4, 5, 6, 7, 8, 1)$, $(1, 2, 3, 4, 5, 6, 7, 8, 9)$, 

$(1, 2, 3, 4, 5, 6, 7, 9, 8)$

The elementary change $\hat 1: 8+ \to 6-$, $\hat 7: 1- \to 1+$ ($8$ enters the conic $23456$, $7$ inside and $1$ outside of $23456$) is realizable only for the
first two pencils. For both, perform the change of pairs: $(1-, 7)$, $(81, C) \to (1+, 7), (1, 6-)$, one gets the two pencils corresponding to the nodal list $\max(\hat 1 = 6-)$.
The elementary change $\hat 1: 6- \to 6+$, $\hat 6: 1- \to 1+$ ($8$ enters the conic $23457$, $1, 6$ outside of $23457$) may be performed on both previous pencils replacing the pair $(6, 1-), (1, 6-)$ by the pair $(6, 1+), (1, 6+)$, the new pencils obtained realize the nodal list $\max(\hat 1 = 6+)$.

\subsection{Proof of Theorem 1}

We will now classify the pencils of cubics with eight base points lying in convex position, up to the action of $D_8$ on these points. Representants of each of the $43$ equivalence classes obtained are shown in Figures~\ref{maxpen}-\ref{pentwo}.
The upper nine pencils in Figure~\ref{maxpen} are the nodal pencils  obtained from the list $\max(\hat 1 = 8-) = (15) \cdot L_{64}$.
The first three of them correspond also to the list  $\max(\hat 1 = 8+) = L_{32}$.
The next two pencils after the blank line correspond to the list  $\max(\hat 1 = 6-) = L_{48}$. The last two pencils after the second blank line correspond to the
list $\max(\hat 1 = 6+) = L_{56}$.
The nodal lists give rise in total to $13$ orbits of pencils, see Figure~\ref{maxpen}.
Recall that two lists that are obtained from one another by an inessential change must be both nodal or both non-nodal. 
For any non-nodal list $L_n$, denote by $\mathcal{P}_n$ the corresponding pencil of cubics.
Hereafter, $(\hat 5, \hat 6)$ stand always for the same elementary change $(\hat 5: 6- \to 6+, \hat 6: 5- \to 5+)$; and $(\hat 1, \hat 8)$ stands for   
 $(\hat 1: 8+ \to 8-, \hat 8: 1+ \to 1-)$.
One has: 
\begin{flushleft}
 $(\hat 5, \hat 6) \cdot L_3 = L_4$,\\ 
$26 \circ (+1)(48) \circ (\hat 1, \hat 8) \cdot L_4 = (+3) \circ (\hat 1, \hat 8) \cdot L_4 = L_{15}$,\\  
$(\hat 5, \hat 6) \circ 37 \cdot L_3 = (\hat 5, \hat 6) \cdot L_{17} = L_{18}$.
\end{flushleft}

Thus $\mathcal{P}_4 = \mathcal{P}_3$;  $\mathcal{P}_{15} = (+3) \cdot \mathcal{P}_4$;  $\mathcal{P}_{18} = 37 \cdot \mathcal{P}_3$.
One has:
\begin{flushleft}
$(\hat 5, \hat 6) \cdot L_5 = L_6$,\\ 
$(\hat 5, \hat 6) \circ 37 \cdot L_5 = (\hat 5, \hat 6) \cdot L_9 = L_{10}$,\\
$26 \circ (+1)(48)  \circ (\hat 1, \hat 8) \cdot L_6 = (+3) \circ (\hat 1, \hat 8) \cdot L_6 = L_{23}$.
\end{flushleft}

Thus,  $\mathcal{P}_6 = \mathcal{P}_5$;  $\mathcal{P}_{10} = 37 \cdot \mathcal{P}_5$;  $\mathcal{P}_{23} = (+3) \cdot \mathcal{P}_6$.
One has $(\hat 5, \hat 6) \circ 37 \cdot L_2 = L_{34}$.  Thus, $\mathcal{P}_{34} = 37 \cdot   \mathcal{P}_2$.
Finally, $(+1)(48) \circ (\hat 1, \hat 8) \cdot L_n = L_m$ for $(n, m) = (7, 26)$, $(8, 25)$, $(11, 22)$, $(12, 21)$, $(13, 20)$, $(14, 19)$. Thus
 $\mathcal{P}_m = (+1)(48) \cdot \mathcal{P}_n$.
The non-nodal principal lists split into two subsets: $28$ lists with $\hat 8 = 1+$ and $15$ lists with $\hat 8 \not= 1+$.
The first set gives rise to $15$ orbits of pencils. In the second set, there is a one-to-one correspondence between the lists and the equivalence classes of pencils,
see Figures~\ref{penone}-\ref{pentwo}, one gets also $15$ pencils.
There are in total $13 + 15 + 15 = 43$ equivalence classes of pencils.

To construct the  non-nodal pencils in the easiest way, we may follow the sequences of elementary changes from Figure~\ref{induct}, with four starting lists:
$L_2$,  $L_{65}$, $L_{71}$, $L_{88}$.
Construct directly the starting pencils corresponding to $L_{65}$, $L_{71}$, $L_{88}$ using again the method exposed in section~3.2.
To get the starting pencil corresponding to $L_2$ in the shortest way,  observe that the list $L_2$ is obtained from the (non-principal and nodal) list $L_1$ by an elementary change $(\hat 6, \hat 7)$.
The list $L_1 = \max(\hat 8 = 1+) = (+1) (48) \max(\hat 1 = 8-)$ corresponds to nine pencils. However, the elementary change $(\hat 6, \hat 7)$ is realizable only for the first of them, shown in the first row of  Figure~\ref{penone}. $\Box$

Pencils of cubics were applied in \cite{F2} to solve an interpolation problem, 
and  in \cite{F1}, \cite{F-O} to study the isotopy types realizable by some real algebraic curves in $\mathbb{R}P^2$. 



\begin{figure}
\begin{tabular}{ c c c c c c c c }
$(12, L)$ & $(81, L)$ & $(81, C)$ & $(1-, 7)$ & $(6, 1-)$ & $(1-, 5)$ & $(4, 1-)$ & $(1-, 3)$\\
$XX$ & $XX$ & $18$ & $X1$ & $61$ & $X1$ & $41$ & $X1$\\
\hline
$(1+, 12)$ & $(1-, 81)$ & $(81, C)$ & $(1-, 7)$ & $(6, 1-)$ & $(1-, 5)$ & $(4, 1-)$ & $(1-, 3)$\\
$8X$ & $2X$ & $X1$ & $18$ & $18$ & $18$ & $18$ & $18$\\
\hline
$(1+, 12)$ & $(1-, 8)$ & $(1-, 78)$ & $(1-, 7)$ & $(6, 1-)$ & $(1-, 5)$ & $(4, 1-)$ & $(1-, 3)$\\
$78$ & $X7$ & $X1$ & $8X$ & $87$ & $87$ & $87$ & $87$\\
\hline
$(1+, 12)$ & $(1-, 8)$ & $(7, 1-)$ & $(1-, 67)$ & $(6, 1-)$ & $(1-, 5)$ & $(4, 1-)$ & $(1-, 3)$\\
$67$ & $76$ & $X6$ & $2X$ & $7X$ & $76$ & $76$ & $76$\\
\hline
$(1+, 12)$ & $(1-, 8)$ & $(7, 1-)$ & $(1-, 6)$ & $(1-, 56)$ & $(1-, 5)$ & $(4, 1-)$ & $(1-, 3)$\\
$56$ & $65$ & $65$ & $X5$ & $X1$ & $6X$ & $65$ & $65$\\
\hline
$(1+, 12)$ & $(1-, 8)$ & $(7, 1-)$ & $(1-, 6)$ & $(5, 1-)$ & $(1-, 45)$ & $(4, 1-)$ & $(1-, 3)$\\
$45$ & $54$ & $54$ & $54$ & $X4$ & $2X$ & $5X$ & $54$\\
\hline
$(1+, 12)$ & $(1-, 8)$ & $(7, 1-)$ & $(1-, 6)$ & $(5, 1-)$ & $(1-, 4)$ & $(1-, 34)$ & $(1-, 3)$\\
$34$ & $43$ & $43$ & $43$ & $43$ & $X3$ & $X1$ & $4X$\\
\hline
$(1+, 12)$ & $(1-, 8)$ & $(7, 1-)$ & $(1-, 6)$ & $(5, 1-)$ & $(1-, 4)$ & $(3, 1-)$ & $(1-, E)$\\
$23$ & $32$ & $32$ & $32$ & $32$ & $32$ & $X2$ & $2X$\\
\hline
$(1+, 12)$ & $(1-, 8)$ & $(7, 1-)$ & $(1-, 6)$ & $(5, 1-)$ & $(1-, 4)$ & $(3, 1-)$ & $(1-, E)$\\
$X2$ & $28$ & $2X$ & $26$ & $2X$ & $24$ & $2X$ & $X2$\\
 & & & & & & & \\
$(12, L)$ & $(81, L)$ & $(1+, 7)$ & $(1, 6-)$ & $(6, 1-)$ & $(1-, 5)$ & $(4, 1-)$ & $(1-, 3)$\\
$XX$ & $XX$ & $X1$ & $16$ & $61$ & $X1$ & $41$ & $X1$\\
\hline
$(1+, 12)$ & $(1-, 18)$ & $(1+, 7)$ & $(1, 6-)$ & $(6, 1-)$ & $(1-, 5)$ & $(4, 1-)$ & $(1-, 3)$\\
$8X$ & $2X$ & $12$ & $X1$ & $18$ & $18$ & $18$ & $18$\\
 & & & & & & & \\
$(12, L)$ & $(81, L)$ & $(1+, 7)$ & $(6, 1+)$ & $(1, 6+)$ & $(1-, 5)$ & $(4, 1-)$ & $(1-, 3)$\\
$XX$ & $XX$ & $X1$ & $61$ & $16$ & $X1$ & $41$ & $X1$\\
\hline
$(1+, 12)$ & $(1-, 18)$ & $(1+, 7)$ & $(6, 1+)$ & $(1, 6+)$ & $(1-, 5)$ & $(4, 1-)$ & $(1-, 3)$\\
$8X$ & $2X$ & $12$ & $12$ & $X1$ & $18$ & $18$ & $18$\\
\end{tabular}
\caption{\label{maxpen} Pencils $\max(\hat 1 = 8-)$, $\max(\hat 1 = 8+)$, $\max(\hat 1 = 6-)$, $\max(\hat 1 = 6+)$}
\end{figure}

\begin{figure}
\begin{tabular}{ c c c c c c c c c }
$(78, L)$ & $(81, L)$ & $(81, C)$ & $(8+, 2)$ & $(3, 8+)$ & $(8+, 4)$ & $(5, 8+)$ & $(8+, 6)$ & $\bf{1}$\\
$XX$ & $XX$ & $18$ & $X8$ & $38$ & $X8$ & $58$ & $X8$ & \\
\hline
$(56, L)$ & $(81, L)$ & $(81, C)$ & $(8+, 2)$ & $(3, 8+)$ & $(8+, 4)$ & $(5, 8+)$ & $(5-, 7)$ & $2$\\
$XX$ & $XX$ & $18$ & $X8$ & $38$ & $X8$ & $58$ & $X5$ & \\
\hline
$(56, L)$ & $(81, L)$ & $(81, C)$ & $(8+, 2)$ & $(3, 8+)$ & $(8+, 4)$ & $(5+, 7)$ & $(56, C)$ & $3, 4$\\
$XX$ & $XX$ & $18$ & $X8$ & $38$ & $X8$ & $X5$ & $65$ & \\
\hline
$(56, L)$ & $(81, L)$ & $(81, C)$ & $(8+, 2)$ & $(3, 8+)$ & $(3-, 7)$ & $(6-, 4)$ & $(56, C)$ & $6, 5$\\
$XX$ & $XX$ & $18$ & $X8$ & $38$ & $X3$ & $X6$ & $65$ & \\
\hline
$(56, L)$ & $(81, L)$ & $(81, C)$ & $(8+, 2)$ & $(3, 8+)$ & $(3-, 7)$ & $(6, 3-)$ & $(6+, 4)$ & $7$\\
$XX$ & $XX$ & $18$ & $X8$ & $38$ & $X3$ & $63$ & $X6$ & \\
\hline
$(34, L)$ & $(81, L)$ & $(81, C)$ & $(8+, 2)$ & $(3, 8+)$ & $(3-, 7)$ & $(6, 3-)$ & $(3-, 5)$ & $8$\\
$XX$ & $XX$ & $18$ & $X8$ & $38$ & $X3$ & $63$ & $X3$ & \\
\hline
$(18, L)$ & $(65, L)$ & $(6+, 4)$ & $(6, 3-)$ & $(3, 6-)$ & $(3+, 7)$ & $(8+, 2)$ & $(18, C)$ & $11$\\
$XX$ & $XX$ & $X6$ & $63$ & $36$ & $X3$ & $X8$ & $81$ & \\
\hline
$(18, L)$ & $(34, L)$ & $(3-, 5)$ & $(6, 3-)$ & $(3, 6-)$ & $(3+, 7)$ & $(8+, 2)$ & $(18, C)$ & $12$\\
$XX$ & $XX$ & $X3$ & $63$ & $36$ & $X3$ & $X8$ & $81$ & \\
\hline
$(18, L)$ & $(34, L)$ & $(3-, 5)$ & $(3, 6+)$ & $(6, 3+)$ & $(3+, 7)$ & $(8+, 2)$ & $(18, C)$ & $14$\\
$XX$ & $XX$ & $X3$ & $36$ & $63$ & $X3$ & $X8$ & $81$ & \\
\hline
$(18, L)$ & $(65, L)$ & $(6+, 4)$ & $(3, 6+)$ & $(6, 3+)$ & $(3+, 7)$ & $(8+, 2)$ & $(18, C)$ & $13$\\
$XX$ & $XX$ & $X6$ & $36$ & $63$ & $X3$ & $X8$ & $81$ & \\
\hline
$(18, L)$ & $(65, L)$ & $(6+, 4)$ & $(6, 3-)$ & $(3, 6-)$ & $(6-, 2)$ & $(1, 6-)$ & $(1+, 7)$ & $35$\\
$XX$ & $XX$ & $X6$ & $63$ & $36$ & $X6$ & $16$ & $X1$\\
\hline
$(18, L)$ & $(34, L)$ & $(3-, 5)$ & $(6, 3-)$ & $(3, 6-)$ & $(6-, 2)$ & $(1, 6-)$ & $(1+, 7)$ & $36$\\
$XX$ & $XX$ & $X3$ & $63$ & $36$ & $X6$ & $16$ & $X1$\\
\hline
$(18, L)$ & $(34, L)$ & $(3-, 5)$ & $(3, 6+)$ & $(6, 3+)$ & $(6-, 2)$ & $(1, 6-)$ & $(1+, 7)$ & $38$\\
$XX$ & $XX$ & $X3$ & $36$ & $63$ & $X6$ & $16$ & $X1$\\
\hline
$(18, L)$ & $(65, L)$ & $(6+, 4)$ & $(3, 6+)$ & $(6, 3+)$ & $(6-, 2)$ & $(1, 6-)$ & $(1+, 7)$ & $37$\\
$XX$ & $XX$ & $X6$ & $36$ & $63$ & $X6$ & $16$ & $X1$\\
\hline
$(18, L)$ & $(65, L)$ & $(6+, 4)$ & $(3, 6+)$ & $(6+, 2)$ & $(6, 1-)$ & $(1, 6-)$ & $(1+, 7)$ & $41$\\
$XX$ & $XX$ & $X6$ & $36$ & $X6$ & $61$ & $16$ & $X1$\\
\hline
$(18, L)$ & $(65, L)$ & $(6+, 4)$ & $(3, 6+)$ & $(6+, 2)$ & $(1, 6+)$ & $(6, 1+)$ & $(1+, 7)$ & $49$\\
$XX$ & $XX$ & $X6$ & $36$ & $X6$ & $16$ & $61$ & $X1$\\
\end{tabular}
\caption{\label{penone} Pencils with $\hat 8 = 1+$}
\end{figure}

\begin{figure}
\begin{tabular}{ c c c c c c c c c }
$(1+, 2)$ & $(3, 1+)$ & $(1+, 4)$ & $(5, 1+)$ & $(1+, 6)$ & $(7, 1+)$ & $(1, 7-)$ & $(1-, 8)$ & $65$\\
$82$ & $8X$ & $84$ & $8X$ & $86$ & $8X$ & $8X$ & $28$\\
\hline
$(1+, 2)$ & $(3, 1+)$ & $(1+, 4)$ & $(5, 1+)$ & $(1+, 6)$ & $(1, 7+)$ & $(7, 1-)$ & $(1-, 8)$ & $66$\\
$82$ & $8X$ & $84$ & $8X$ & $86$ & $6X$ & $2X$ & $28$\\
\hline
$(1+, 2)$ & $(3, 1+)$ & $(1+, 4)$ & $(5, 1+)$ & $(1, 5-)$ & $(1-, 6)$ & $(7, 1-)$ & $(1-, 8)$ & $67$\\
$82$ & $8X$ & $84$ & $8X$ & $6X$ & $26$ & $2X$ & $28$\\
\hline
$(4, 2+)$ & $(2+, 5)$ & $(6, 2+)$ & $(2+, 7)$ & $(8, 2+)$ & $(2, 8-)$ & $(8-, 3)$ & $(4+, 1)$ & $71$\\
$1X$ & $15$ & $1X$ & $17$ & $1X$ & $1X$ & $13$ & $31$\\
\hline
$(4-, 1)$ & $(2+, 5)$ & $(6, 2+)$ & $(2+, 7)$ & $(8, 2+)$ & $(2, 8-)$ & $(8-, 3)$ & $(4, 8-)$ & $72$\\
$51$ & $15$ & $1X$ & $17$ & $1X$ & $1X$ & $13$ & $1X$\\
\hline
$(4, 2+)$ & $(2+, 5)$ & $(6, 2+)$ & $(2+, 7)$ & $(2, 8+)$ & $(8, 2-)$ & $(8-, 3)$ & $(4+, 1)$ & $75$\\
$1X$ & $15$ & $1X$ & $17$ & $7X$ & $3X$ & $13$ & $31$\\
\hline
$(4, 2+)$ & $(2+, 5)$ & $(6, 2+)$ & $(2, 6-)$ & $(2-, 7)$ & $(8, 2-)$ & $(8-, 3)$ & $(4+, 1)$ & $78$\\
$1X$ & $15$ & $1X$ & $7X$ & $37$ & $3X$ & $13$ & $31$\\
\hline
$(4-, 1)$ & $(2+, 5)$ & $(6, 2+)$ & $(2+, 7)$ & $(2, 8+)$ & $(8, 2-)$ & $(8-, 3)$ & $(4, 8-)$ & $80$\\
$51$ & $15$ & $1X$ & $17$ & $7X$ & $3X$ & $13$ & $1X$\\
\hline
$(4-, 1)$ & $(2+, 5)$ & $(6, 2+)$ & $(2, 6-)$ & $(2-, 7)$ & $(8, 2-)$ & $(8-, 3)$ & $(4, 8-)$ & $86$\\
$51$ & $15$ & $1X$ & $7X$ & $37$ & $3X$ & $13$ & $1X$\\
\hline
$(4-, 1)$ & $(2+, 5)$ & $(2, 6+)$ & $(6, 2-)$ & $(2-, 7)$ & $(8, 2-)$ & $(8-, 3)$ & $(4, 8-)$ & $87$\\
$51$ & $15$ & $5X$ & $3X$ & $37$ & $3X$ & $13$ & $1X$\\
\hline
$(4-, 1)$ & $(2, 4-)$ & $(2-, 5)$ & $(6, 2-)$ & $(2-, 7)$ & $(8, 2-)$ & $(8-, 3)$ & $(4, 8-)$ & $84$\\
$51$ & $5X$ & $35$ & $3X$ & $37$ & $3X$ & $13$ & $1X$\\
\hline
$(8-, 5)$ & $(6+, 1)$ & $(2, 6+)$ & $(6, 2-)$ & $(2-, 7)$ & $(8, 2-)$ & $(8-, 3)$ & $(4, 8-)$ & $83$\\
$15$ & $51$ & $5X$ & $3X$ & $37$ & $3X$ & $13$ & $1X$\\
\hline
$(8-, 5)$ & $(6+, 1)$ & $(6, 2+)$ & $(2, 6-)$ & $(2-, 7)$ & $(8, 2-)$ & $(8-, 3)$ & $(4, 8-)$ & $82$\\
$15$ & $51$ & $1X$ & $7X$ & $37$ & $3X$ & $13$ & $1X$\\
\hline
$(5, 1-)$ & $(1, 5+)$ & $(1+, 4)$ & $(3-, 8)$ & $(3, 7-)$ & $(7, 3+)$ & $(7+, 2)$ & $(1-, 6)$ & $88$\\
$2X$ & $4X$ & $84$ & $48$ & $8X$ & $2X$ & $62$ & $26$\\
\hline
$(1, 5-)$ & $(5, 1+)$ & $(1+, 4)$ & $(3-, 8)$ & $(7, 3-)$ & $(3, 7+)$ & $(7+, 2)$ & $(1-, 6)$ & $92$\\
$6X$ & $8X$ & $84$ & $48$ & $4X$ & $6X$ & $62$ & $26$\\
\end{tabular}
\caption{\label{pentwo} Pencils with $\hat 8 \not= 1+$}
\end{figure}

\newpage
\section{Tabulars}


\begin{figure}
\begin{tabular}{|c|c|c|c|c|c|c|c|c|c|c|c|c|c|c|}
\hline
$\hat 8$ & \multicolumn{2}{c|}{$1+$} & \multicolumn{2}{c|}{$1-$}  &  \multicolumn{2}{c|}{$2+$} &  \multicolumn{2}{c|}{$2-$}
&  \multicolumn{2}{c|}{$3+$} &  \multicolumn{2}{c|}{$3-$} & \multicolumn{2}{c|}{$4+$}\\
\hline
$C_2$ & in & out & in & out & in & out & in & out & in & out & in & out & in & out\\
$23456$ & ${\bf 7}$ & ${\bf 1}$ &  & $1, 7$ & $1, 7$ &  & $1, 7$ &  & & $1, 7$ & & $1, 7$ & $1, 7$ & \\
$23457$ &  & $1, 6$ & $6$ & $1$ & $1$ & $6$ & $1$ & $6$ & $6$ & $1$ & $6$ & $1$ & $1$ & $6$\\
$23467$ & $5$ & $1$ &  & $1, 5$ & $1, 5$ &   & $1, 5$ & & & $1, 5$ & & $1, 5$ & $1, 5$ & \\
$23567$ &  & $1, 4$ & $4$ & $1$ & $1$ & $4$ & $1$ & $4$ & $4$ & $1$ & $4$ & $1$ & $1$ & $4$\\
$24567$ & $3$ & $1$ &  & $1, 3$ & $1, 3$ &  & $1, 3$ & & & $1, 3$ & $1$ & $3$ & $3$ & $1$\\
$34567$ &  & $1, 2$ & ${\bf 2}$ & ${\bf 1}$ & ${\bf 1}$ & ${\bf 2}$ &  & $1, 2$ & $1, 2$ & & $1, 2$ & & & $1, 2$\\
$13456$ & $2, 7$ &  & $2, 7$ &  &  & $2, 7$ & $7$ & $2$ & $2$ & $7$ & $2$ & $7$ & $7$ & $2$\\
$13457$ & $2$ & $6$ & $2$ & $6$ & $6$ & $2$ &  & $2, 6$ & $2, 6$ & & $2, 6$ & & & $2, 6$\\
$13467$ & $2, 5$ &  & $2, 5$ &  & & $2, 5$ & $5$ & $2$ & $2$ & $5$ & $2$ & $5$ & $5$ & $2$\\
$13567$ & $2$ & $4$ & $2$ & $4$ & $4$ & $2$ & & $2, 4$ & $2, 4$ & & $2, 4$ & & & $2, 4$\\
$14567$ & $2, 3$ &  & $2, 3$ &  & & $2, 3$ & ${\bf 3}$ & ${\bf 2}$ & ${\bf 2}$ & ${\bf 3}$ & & $2, 3$ & $2, 3$ & \\
$12456$ & $7$ & $3$ & $7$ & $3$ & $3$ & $7$ & $3$ & $7$ & $7$ & $3$ & & $3, 7$ & $3, 7$ & \\
$12457$ &  & $3, 6$ &  & $3, 6$ & $3, 6$ & & $3, 6$ & & & $3, 6$ & $6$ & $3$ & $3$ & $6$\\
$12467$ & $5$ & $3$ & $5$ & $3$ & $3$ & $5$ & $3$ & $5$ & $5$ &  $3$ & & $3, 5$ & $3, 5$ & \\
$12567$ &  & $3, 4$ &  & $3, 4$ & $3, 4$ & & $3, 4$ & & & $3, 4$ & ${\bf 4}$ & ${\bf 3}$ & ${\bf 3}$ & ${\bf 4}$\\
$12356$ & $4, 7$ &  & $4, 7$ &  & & $4, 7$ & & $4, 7$ & $4, 7$ & & $4, 7$ & & & $4, 7$\\
$12357$ & $4$ & $6$ & $4$ & $6$ & $6$ & $4$ & $6$ & $4$ & $4$ & $6$ & $4$ & $6$ & $6$ & $4$\\
$12367$ & $4, 5$ &  & $4, 5$ & & & $4, 5$ & & $4, 5$ & $4, 5$ & & $4, 5$ & & & $4, 5$\\
$12346$ & $7$ & $5$ & $7$ & $5$ & $5$ & $7$ & $5$ & $7$ & $7$ & $5$ & $7$ & $5$ & $5$ & $7$\\
$12347$ &  & $5, 6$ &  & $5, 6$ & $5, 6$ & & $5, 6$ & & & $5, 6$ & & $5, 6$ & $5, 6$ & \\
$12345$ & $6, 7$ &  & $6, 7$ & & & $6, 7$ & & $6, 7$ & $6, 7$ & & $6, 7$ & & & $6, 7$\\
\hline
\end{tabular}
\caption{\label{lsevenone} The lists $\hat 8 = L(1, \dots 7)$}
\end{figure}

\begin{figure}
\begin{tabular}{|c|c|c|c|c|c|c|c|c|c|c|c|c|c|c|}
\hline
$\hat 8$ & \multicolumn{2}{c|}{$4-$} & \multicolumn{2}{c|}{$5+$}  &  \multicolumn{2}{c|}{$5-$} &  \multicolumn{2}{c|}{$6+$}
&  \multicolumn{2}{c|}{$6-$} &  \multicolumn{2}{c|}{$7+$} & \multicolumn{2}{c|}{$7-$}\\
\hline
$C_2$ & in & out & in & out & in & out & in & out & in & out & in & out & in & out\\
$23456$ & $1, 7$ & & & $1, 7$ & & $1, 7$ & $1, 7$ & & $1, 7$ & & & $1, 7$ & ${\bf 1}$ & ${\bf 7}$\\
$23457$ & $1$ & $6$ & $6$ & $1$ & $6$ & $1$ & $1$ & $6$ & & $1, 6$ & $1, 6$ & & $1, 6$ & \\
$23467$ & $1, 5$ & & & $1, 5$ & $1$ & $5$ & $5$ & $1$ & $5$ & $1$ & $1$ & $5$ & $1$ & $5$\\
$23567$ & & $1, 4$ & $1, 4$ & & $1, 4$ & & & $1, 4$ & & $1, 4$ & $1, 4$ & & $1, 4$ & \\
$24567$ & $3$ & $1$ & $1$ & $3$ & $1$ & $3$ & $3$ & $1$ & $3$ & $1$ & $1$ & $3$ & $1$ & $3$\\
$34567$ & & $1, 2$ & $1, 2$ & & $1, 2$ & & & $1, 2$ & & $1, 2$ & $1, 2$ & & $1, 2$ & \\
$13456$ & $7$ & $2$ & $2$ & $7$ & $2$ & $7$ & $7$ & $2$ & $7$ & $2$ & $2$ & $7$ & & $2, 7$\\
$13457$ & & $2, 6$ & $2, 6$ & & $2, 6$ & & & $2, 6$ & $2$ & $6$ & $6$ & $2$ & $6$ & $2$\\
$13467$ & $5$ & $2$ & $2$ & $5$ & & $2, 5$ & $2, 5$ & & $2, 5$ & & & $2, 5$ & & $2, 5$\\
$13567$ & $2$ & $4$ & $4$ & $2$ & $4$ & $2$ & $2$ & $4$ & $2$ & $4$ & $4$ & $2$ & $4$ & $2$\\
$14567$ & $2, 3$ & & & $2, 3$ & & $2, 3$ & $2, 3$ & & $2, 3$ & & & $2, 3$ & & $2, 3$\\
$12456$ & $3, 7$ & & & $3, 7$ & & $3, 7$ & $3, 7$ & & $3, 7$ & & & $3, 7$ & $3$ & $7$\\
$12457$ & $3$ & $6$ & $6$ & $3$ & $6$ & $3$ & $3$ & $6$ & & $3, 6$ & $3, 6$ & & $3, 6$ & \\
$12467$ & $3, 5$ & & & $3, 5$ & $3$ & $5$ & $5$ & $3$ & $5$ & $3$ & $3$ & $5$ & $3$ & $5$\\
$12567$ & & $3, 4$ & $3, 4$ & & $3, 4$ & & & $3, 4$ & & $3, 4$ & $3, 4$ & & $3, 4$ & \\
$12356$ & $7$ & $4$ & $4$ & $7$ & $4$ & $7$ & $7$ & $4$ & $7$ & $4$ & $4$ & $7$ & & $4, 7$\\
$12357$ & & $4, 6$ & $4, 6$ & & $4, 6$ & & & $4, 6$ & $4$ & $6$ & $6$ & $4$ & $6$ & $4$\\
$12367$ & ${\bf 5}$ & ${\bf 4}$ & ${\bf 4}$ & ${\bf 5}$ & & $4, 5$ & $4, 5$ & & $4, 5$ & & & $4, 5$ & & $4, 5$\\
$12346$ & $5$ & $7$ & $7$ & $5$ & & $5, 7$ & $5, 7$ & & $5, 7$ & & & $5, 7$ & $5$ & $7$\\
$12347$ & $5, 6$ & & & $5, 6$ & ${\bf 6}$ & ${\bf 5}$ & ${\bf 5}$ & ${\bf 6}$ & & $5, 6$ & $5, 6$ & & $5, 6$ & \\
$12345$ & & $6, 7$ & $6, 7$ & & $6, 7$ & & & $6, 7$ & ${\bf 7}$ & ${\bf 6}$ & ${\bf 6}$ & ${\bf 7}$ & & $6, 7$\\
\hline
\end{tabular}
\caption{\label{lseventwo} The lists $\hat 8 = L(1, \dots 7)$, continued}
\end{figure}

\begin{figure}
\begin{displaymath}
\xymatrix{
\hat 1  \ar[r]^{8+} &  23456 \ar[r]^{6-} & 23457  \ar[r]^{6+} & 23467 \ar[r]^{4-} & 23567 \ar[r]^{4+} & 24567 
 \ar[r]^{2-} & 34567 \ar[r]^{2+} \ar@{=}[d] & \hat 1\\
\hat 2  \ar[r]^{8+} &  13456 \ar[r]^{6-}  \ar[u] & 13457  \ar[r]^{6+} \ar[u]  & 13467 \ar[r]^{4-} \ar[u] 
& 13567 \ar[r]^{4+} \ar[u]  & 14567  \ar[r]^{1-}  \ar[u] \ar@{=}[d] & 34567 \ar[r]^{1+} & \hat 2\\
\hat 3  \ar[r]^{8+} &  12456 \ar[r]^{6-}  \ar[u] & 12457  \ar[r]^{6+} \ar[u]  & 12467 \ar[r]^{4-} \ar[u]
& 12567 \ar[r]^{4+} \ar[u]  \ar@{=}[d] & 14567  \ar[r]^{1-} & 24567 \ar[r]^{1+}  \ar[u]  & \hat 3\\
\hat 4  \ar[r]^{8+} &  12356 \ar[r]^{6-}  \ar[u] & 12357  \ar[r]^{6+} \ar[u] & 12367 \ar[r]^{3-} \ar[u] \ar@{=}[d]
& 12567 \ar[r]^{3+} & 13567 \ar[r]^{1-}  \ar[u]   & 23567 \ar[r]^{1+}  \ar[u] & \hat 4\\
\hat 5  \ar[r]^{8+} &  12346 \ar[r]^{6-} \ar[u] & 12347  \ar[r]^{6+} \ar[u]  \ar@{=}[d] & 12367 \ar[r]^{3-} & 12467 \ar[r]^{3+} 
\ar[u]   & 13467 \ar[r]^{1-} \ar[u] & 23467 \ar[r]^{1+}  \ar[u]  & \hat 5\\
\hat 6  \ar[r]^{8+} &  12345 \ar[r]^{5-} \ar[u]  \ar@{=}[d] & 12347  \ar[r]^{5+} & 12357 \ar[r]^{3-} \ar[u]   & 12457 \ar[r]^{3+} 
\ar[u]  & 13457 \ar[r]^{1-}  \ar[u]  & 23457 \ar[r]^{1+}  \ar[u]   & \hat 6\\
\hat 7  \ar[r]^{8+} &  12345 \ar[r]^{5-} & 12346  \ar[r]^{5+} \ar[u] & 12356 \ar[r]^{3-}  \ar[u]  & 12456 \ar[r]^{3+} 
\ar[u]   & 13456 \ar[r]^{1-}  \ar[u]  & 23456 \ar[r]^{1+} \ar[u] & \hat 7\\
}
\end{displaymath}
\caption{\label{oneplus} $\hat 8 = 1+$}
\end{figure}

\begin{figure}
\begin{displaymath}
\xymatrix{
\hat 1  \ar[r]^{8-} &  34567 \ar[r]^{3+} & 24567  \ar[r]^{3-} & 23567 \ar[r]^{5+} & 23467 \ar[r]^{5-} & 23457 
 \ar[r]^{7+} & 23456 \ar[r]^{7-} & \hat 1\\
\hat 2  \ar[r]^{8+} &  13456 \ar[r]^{6-}  \ar[u] & 13457  \ar[r]^{6+} \ar[u]  & 13467 \ar[r]^{4-} \ar[u] 
& 13567 \ar[r]^{4+} \ar[u]  & 14567  \ar[r]^{1-}  \ar[u] \ar@{=}[d] & 34567 \ar[r]^{1+}  \ar[u] \ar[d] & \hat 2\\
\hat 3  \ar[r]^{8+} &  12456 \ar[r]^{6-}  \ar[u] & 12457  \ar[r]^{6+} \ar[u]  & 12467 \ar[r]^{4-} \ar[u]
& 12567 \ar[r]^{4+} \ar[u] \ar@{=}[d] & 14567  \ar[r]^{1-} & 24567 \ar[r]^{1+}  \ar[d]  & \hat 3\\
\hat 4  \ar[r]^{8+} &  12356 \ar[r]^{6-}  \ar[u] & 12357  \ar[r]^{6+} \ar[u] & 12367 \ar[r]^{3-} \ar[u] \ar@{=}[d]
& 12567 \ar[r]^{3+} & 13567 \ar[r]^{1-}  \ar[u]   & 23567 \ar[r]^{1+}  \ar[d] & \hat 4\\
\hat 5 \ar[r]^{8+} &  12346 \ar[r]^{6-} \ar[u] & 12347  \ar[r]^{6+} \ar[u]  \ar@{=}[d] & 12367 \ar[r]^{3-} & 12467 \ar[r]^{3+} 
\ar[u]   & 13467 \ar[r]^{1-} \ar[u] & 23467 \ar[r]^{1+}  \ar[d]  & \hat 5\\
\hat 6 \ar[r]^{8+} &  12345 \ar[r]^{5-} \ar[u] \ar@{=}[d] & 12347  \ar[r]^{5+} & 12357 \ar[r]^{3-} \ar[u]   & 12457 \ar[r]^{3+} 
\ar[u]  & 13457 \ar[r]^{1-}  \ar[u]  & 23457 \ar[r]^{1+}  \ar[d]   & \hat 6\\
\hat 7 \ar[r]^{8+} &  12345 \ar[r]^{5-} & 12346  \ar[r]^{5+} \ar[u] & 12356 \ar[r]^{3-}  \ar[u]  & 12456 \ar[r]^{3+} 
\ar[u]   & 13456 \ar[r]^{1-}  \ar[u]  & 23456 \ar[r]^{1+} & \hat 7\\
}
\end{displaymath}
\caption{\label{oneminus} $\hat 8 = 1-$}
\end{figure}

\begin{figure}
\begin{displaymath}
\xymatrix{
\hat 1  \ar[r]^{8+} &  23456 \ar[r]^{6-}  \ar[d] & 23457  \ar[r]^{6+}  \ar[d] & 23467 \ar[r]^{4-} \ar[d] & 23567 \ar[r]^{4+} \ar[d] & 
24567  \ar[r]^{2-} \ar[d] & 34567 \ar[r]^{2+} \ar[d] & \hat 1\\
\hat 2  \ar[r]^{8-} &  34567 \ar[r]^{3+} & 14567  \ar[r]^{3-} \ar@{=}[d]  & 13567 \ar[r]^{5+} \ar[d] 
& 13467 \ar[r]^{5-} \ar[d]  & 13457  \ar[r]^{7+}  \ar[d] & 13456 \ar[r]^{7-}  \ar[d] & \hat 2\\
\hat 3  \ar[r]^{8-} &  24567 \ar[r]^{2+}  \ar[u] & 14567  \ar[r]^{2-} \ar[d]  & 12567 \ar[r]^{5+} \ar@{=}[d]
& 12467 \ar[r]^{5-} \ar[d] & 12457  \ar[r]^{7+} \ar[d]  & 12456 \ar[r]^{7-}  \ar[d]  & \hat 3\\
\hat 4  \ar[r]^{8-}  &  23567 \ar[r]^{2+}  \ar[u] & 13567  \ar[r]^{2-} \ar[d] & 12567 \ar[r]^{5+} \ar[d] 
& 12367 \ar[r]^{5-} \ar@{=}[d] & 12357 \ar[r]^{7+}  \ar[d]   & 12356 \ar[r]^{7-}  \ar[d] & \hat 4\\
\hat 5 \ar[r]^{8-} &  23467 \ar[r]^{2+} \ar[u] & 13467  \ar[r]^{2-} \ar[d]  & 12467 \ar[r]^{4+} \ar[d] & 12367 \ar[r]^{4-} 
\ar[d]  & 12347 \ar[r]^{7+} \ar@{=}[d] & 12346 \ar[r]^{7-}  \ar[d]  & \hat 5\\
\hat 6 \ar[r]^{8-} &  23457 \ar[r]^{2+} \ar[u] & 13457  \ar[r]^{2-} \ar[d]  & 12457 \ar[r]^{4+} \ar[d]  & 12357 \ar[r]^{4-} 
\ar[d]  & 12347 \ar[r]^{7+}  \ar[d]  & 12345 \ar[r]^{7-}  \ar@{=}[d]  & \hat 6\\
\hat 7 \ar[r]^{8-} &  23456 \ar[r]^{2+} \ar[u] & 13456  \ar[r]^{2-} & 12456 \ar[r]^{4+} & 12356 \ar[r]^{4-} 
& 12346 \ar[r]^{6+}  & 12345 \ar[r]^{6-} & \hat 7\\
}
\end{displaymath}
\caption{\label{twoplus} $\hat 8 = 2+$}
\end{figure}

\begin{figure}
\begin{displaymath}
\xymatrix{
\hat 1  \ar[r]^{8+} &  23456 \ar[r]^{6-}  \ar[d] & 23457  \ar[r]^{6+}  \ar[d] & 23467 \ar[r]^{4-} \ar[d] & 23567 \ar[r]^{4+} \ar[d] & 
24567  \ar[r]^{2-} \ar[d] & 34567 \ar[r]^{2+} \ar@{=}[d] & \hat 1\\
\hat 2  \ar[r]^{8+} &  13456 \ar[r]^{6-}  & 13457  \ar[r]^{6+} \ar[d]  & 13467 \ar[r]^{4-} \ar[d] 
& 13567 \ar[r]^{4+} \ar[d]  & 14567  \ar[r]^{1-}  \ar[d] & 34567 \ar[r]^{1+}  & \hat 2\\
\hat 3  \ar[r]^{8-} &  24567 \ar[r]^{2+}  & 14567  \ar[r]^{2-} & 12567 \ar[r]^{5+} \ar@{=}[d]
& 12467 \ar[r]^{5-} \ar[d] & 12457  \ar[r]^{7+} \ar[d]  & 12456 \ar[r]^{7-}  \ar[d]  & \hat 3\\
\hat 4  \ar[r]^{8-}  &  23567 \ar[r]^{2+}  \ar[u] & 13567  \ar[r]^{2-} \ar[u] & 12567 \ar[r]^{5+} \ar[d] 
& 12367 \ar[r]^{5-} \ar@{=}[d] & 12357 \ar[r]^{7+}  \ar[d]   & 12356 \ar[r]^{7-}  \ar[d] & \hat 4\\
\hat 5 \ar[r]^{8-} &  23467 \ar[r]^{2+} \ar[u] & 13467  \ar[r]^{2-} \ar[u]  & 12467 \ar[r]^{4+} \ar[d] & 12367 \ar[r]^{4-} 
\ar[d]  & 12347 \ar[r]^{7+} \ar@{=}[d] & 12346 \ar[r]^{7-}  \ar[d]  & \hat 5\\
\hat 6 \ar[r]^{8-} &  23457 \ar[r]^{2+} \ar[u] & 13457  \ar[r]^{2-} \ar[u]  & 12457 \ar[r]^{4+} \ar[d]  & 12357 \ar[r]^{4-} 
\ar[d]  & 12347 \ar[r]^{7+}  \ar[d]  & 12345 \ar[r]^{7-}  \ar@{=}[d]  & \hat 6\\
\hat 7 \ar[r]^{8-} &  23456 \ar[r]^{2+} \ar[u] & 13456  \ar[r]^{2-} \ar[u] & 12456 \ar[r]^{4+} & 12356 \ar[r]^{4-} 
& 12346 \ar[r]^{6+}  & 12345 \ar[r]^{6-} & \hat 7\\
}
\end{displaymath}
\caption{\label{twominus} $\hat 8 = 2-$}
\end{figure}

\begin{figure}
\begin{displaymath}
\xymatrix{
\hat 1  \ar[r]^{8-} &  34567 \ar[r]^{3+} \ar@{=}[d] & 24567  \ar[r]^{3-} & 23567 \ar[r]^{5+} & 23467 \ar[r]^{5-} & 23457 
 \ar[r]^{7+} & 23456 \ar[r]^{7-} & \hat 1\\
\hat 2  \ar[r]^{8-} &  34567 \ar[r]^{3+} & 14567  \ar[r]^{3-} \ar[u]  & 13567 \ar[r]^{5+} \ar[u] 
& 13467 \ar[r]^{5-} \ar[u]  & 13457  \ar[r]^{7+}  \ar[u] & 13456 \ar[r]^{7-}  \ar[u] & \hat 2\\
\hat 3  \ar[r]^{8+} &  12456 \ar[r]^{6-}  & 12457  \ar[r]^{6+} \ar[u]  & 12467 \ar[r]^{4-} \ar[u]
& 12567 \ar[r]^{4+} \ar[u] \ar@{=}[d] & 14567  \ar[r]^{1-} \ar[u] \ar[d] & 24567 \ar[r]^{1+}  \ar[d]  & \hat 3\\
\hat 4  \ar[r]^{8+} &  12356 \ar[r]^{6-}  \ar[u] & 12357  \ar[r]^{6+} \ar[u] & 12367 \ar[r]^{3-} \ar[u] \ar@{=}[d]
& 12567 \ar[r]^{3+} & 13567 \ar[r]^{1-}  \ar[d]   & 23567 \ar[r]^{1+}  \ar[d] & \hat 4\\
\hat 5 \ar[r]^{8+} &  12346 \ar[r]^{6-} \ar[u] & 12347  \ar[r]^{6+} \ar[u]  \ar@{=}[d] & 12367 \ar[r]^{3-} & 12467 \ar[r]^{3+} 
\ar[u]   & 13467 \ar[r]^{1-} \ar[d] & 23467 \ar[r]^{1+}  \ar[d]  & \hat 5\\
\hat 6 \ar[r]^{8+} &  12345 \ar[r]^{5-} \ar[u] \ar@{=}[d] & 12347  \ar[r]^{5+} & 12357 \ar[r]^{3-} \ar[u]   & 12457 \ar[r]^{3+} 
\ar[u]  & 13457 \ar[r]^{1-}  \ar[d]  & 23457 \ar[r]^{1+}  \ar[d]   & \hat 6\\
\hat 7 \ar[r]^{8+} &  12345 \ar[r]^{5-} & 12346  \ar[r]^{5+} \ar[u] & 12356 \ar[r]^{3-}  \ar[u]  & 12456 \ar[r]^{3+} 
\ar[u]   & 13456 \ar[r]^{1-}  & 23456 \ar[r]^{1+} & \hat 7\\
}
\end{displaymath}
\caption{\label{threeplus} $\hat 8 = 3+$}
\end{figure}

\begin{figure}
\begin{displaymath}
\xymatrix{
\hat 1  \ar[r]^{8-} &  34567 \ar[r]^{3+} \ar@{=}[d] & 24567  \ar[r]^{3-} \ar[d] & 23567 \ar[r]^{5+} & 23467 \ar[r]^{5-} & 23457 
 \ar[r]^{7+} & 23456 \ar[r]^{7-} & \hat 1\\
\hat 2  \ar[r]^{8-} &  34567 \ar[r]^{3+} \ar[d] & 14567  \ar[r]^{3-} \ar@{=}[d]  & 13567 \ar[r]^{5+} \ar[u] 
& 13467 \ar[r]^{5-} \ar[u]  & 13457  \ar[r]^{7+}  \ar[u] & 13456 \ar[r]^{7-}  \ar[u] & \hat 2\\
\hat 3  \ar[r]^{8-} &  24567 \ar[r]^{2+}  & 14567  \ar[r]^{2-} & 12567 \ar[r]^{5+} \ar[u]
& 12467 \ar[r]^{5-} \ar[u] & 12457  \ar[r]^{7+} \ar[u]  & 12456 \ar[r]^{7-}  \ar[u]  & \hat 3\\
\hat 4  \ar[r]^{8+} &  12356 \ar[r]^{6-} & 12357  \ar[r]^{6+} & 12367 \ar[r]^{3-} \ar[u] \ar@{=}[d]
& 12567 \ar[r]^{3+} \ar[u] \ar[d] & 13567 \ar[r]^{1-}  \ar[d]   & 23567 \ar[r]^{1+}  \ar[d] & \hat 4\\
\hat 5 \ar[r]^{8+} &  12346 \ar[r]^{6-} \ar[u] & 12347  \ar[r]^{6+} \ar[u]  \ar@{=}[d] & 12367 \ar[r]^{3-} & 12467 \ar[r]^{3+} 
\ar[d]   & 13467 \ar[r]^{1-} \ar[d] & 23467 \ar[r]^{1+}  \ar[d]  & \hat 5\\
\hat 6 \ar[r]^{8+} &  12345 \ar[r]^{5-} \ar[u] \ar@{=}[d] & 12347  \ar[r]^{5+} & 12357 \ar[r]^{3-} \ar[u] & 12457 \ar[r]^{3+} 
\ar[d]  & 13457 \ar[r]^{1-}  \ar[d]  & 23457 \ar[r]^{1+}  \ar[d]   & \hat 6\\
\hat 7 \ar[r]^{8+} &  12345 \ar[r]^{5-} & 12346  \ar[r]^{5+} \ar[u] & 12356 \ar[r]^{3-}  \ar[u]  & 12456 \ar[r]^{3+} 
& 13456 \ar[r]^{1-}  & 23456 \ar[r]^{1+} & \hat 7\\
}
\end{displaymath}
\caption{\label{threeminus} $\hat 8 = 3-$}
\end{figure}

\begin{figure}
\begin{displaymath}
\xymatrix{
\hat 1  \ar[r]^{8+} &  23456 \ar[r]^{6-}  \ar[d] & 23457  \ar[r]^{6+} \ar[d] & 23467 \ar[r]^{4-} \ar[d] & 23567 \ar[r]^{4+} \ar[d] & 
24567  \ar[r]^{2-} & 34567 \ar[r]^{2+} \ar@{=}[d] & \hat 1\\
\hat 2  \ar[r]^{8+} &  13456 \ar[r]^{6-} \ar[d] & 13457  \ar[r]^{6+} \ar[d]  & 13467 \ar[r]^{4-} \ar[d] 
& 13567 \ar[r]^{4+} \ar[d]  & 14567  \ar[r]^{1-} \ar[u] \ar@{=}[d] & 34567 \ar[r]^{1+}  & \hat 2\\
\hat 3  \ar[r]^{8+} &  12456 \ar[r]^{6-}  & 12457  \ar[r]^{6+} & 12467 \ar[r]^{4-} \ar[d]
& 12567 \ar[r]^{4+} \ar[d] & 14567  \ar[r]^{1-} & 24567 \ar[r]^{1+}  \ar[u]  & \hat 3\\
\hat 4  \ar[r]^{8-}  &  23567 \ar[r]^{2+}  & 13567  \ar[r]^{2-} & 12567 \ar[r]^{5+}  
& 12367 \ar[r]^{5-} \ar@{=}[d] & 12357 \ar[r]^{7+}  \ar[d]   & 12356 \ar[r]^{7-}  \ar[d] & \hat 4\\
\hat 5 \ar[r]^{8-} &  23467 \ar[r]^{2+} \ar[u] & 13467  \ar[r]^{2-} \ar[u]  & 12467 \ar[r]^{4+} \ar[u] & 12367 \ar[r]^{4-} 
\ar[d]  & 12347 \ar[r]^{7+} \ar@{=}[d] & 12346 \ar[r]^{7-}  \ar[d]  & \hat 5\\
\hat 6 \ar[r]^{8-} &  23457 \ar[r]^{2+} \ar[u] & 13457  \ar[r]^{2-} \ar[u]  & 12457 \ar[r]^{4+} \ar[u]  & 12357 \ar[r]^{4-} 
\ar[d]  & 12347 \ar[r]^{7+}  \ar[d]  & 12345 \ar[r]^{7-}  \ar@{=}[d]  & \hat 6\\
\hat 7 \ar[r]^{8-} &  23456 \ar[r]^{2+} \ar[u] & 13456  \ar[r]^{2-} \ar[u] & 12456 \ar[r]^{4+} \ar[u] & 12356 \ar[r]^{4-} 
& 12346 \ar[r]^{6+}  & 12345 \ar[r]^{6-} & \hat 7\\
}
\end{displaymath}
\caption{\label{fourplus} $\hat 8 = 4+$}
\end{figure}


\begin{figure}
\begin{tabular}{c c c c c}
$\hat 1$ & $3+ \leftrightarrow 8-$ & $3+ \leftrightarrow   8-$ & $2+ \leftrightarrow 2-$ & $2+ \leftrightarrow 2-$\\ 
$\hat 2$ & $3+ \leftrightarrow 8-$ & $1+ \leftrightarrow   1-$ & $3+ \leftrightarrow 8-$ & $1+ \leftrightarrow 1-$\\ 
 & & \\
$\hat 1$ & $3- \leftrightarrow 3+$ & $3- \leftrightarrow 3+$ & $2- \leftrightarrow 4+$ & $2- \leftrightarrow 4+$\\ 
$\hat 3$ & $2+ \leftrightarrow 8-$ & $1+ \leftrightarrow 1-$ & $2+ \leftrightarrow 8-$ & $1+ \leftrightarrow 1-$\\ 
 & & \\
$\hat 1$ & $5+  \leftrightarrow 3-$ & $5+  \leftrightarrow 3-$ & $4+  \leftrightarrow 4-$ & $4+  \leftrightarrow 4-$\\ 
$\hat 4$ & $2+  \leftrightarrow 8-$ & $1+  \leftrightarrow 1-$ & $2+  \leftrightarrow 8-$ & $1+  \leftrightarrow 1-$\\ 
 & & \\
$\hat 1$ & $5-  \leftrightarrow 5+$ & $5-  \leftrightarrow 5+$ & $4-  \leftrightarrow 6+$ & $4-  \leftrightarrow 6+$\\ 
$\hat 5$ & $2+  \leftrightarrow 8-$ & $1+  \leftrightarrow 1-$ & $2+  \leftrightarrow 8-$ & $1+  \leftrightarrow 1-$\\ 
 & & \\
$\hat 1$ & $7+  \leftrightarrow 5-$ & $7+  \leftrightarrow 5-$ & $6+  \leftrightarrow 6-$ & $6+  \leftrightarrow 6-$\\ 
$\hat 6$ & $2+  \leftrightarrow 8-$ & $1+  \leftrightarrow 1-$ & $2+  \leftrightarrow 8-$ & $1+  \leftrightarrow 1-$\\ 
 & & \\
$\hat 1$ & $7-  \leftrightarrow 7+$ & $7-  \leftrightarrow 7+$ & $6-  \leftrightarrow 8+$ & $6-  \leftrightarrow 8+$\\ 
$\hat 7$ & $2+  \leftrightarrow 8-$ & $1+  \leftrightarrow 1-$ & $2+  \leftrightarrow 8-$ & $1+  \leftrightarrow 1-$\\ 
 & & \\
$\hat 1$ & $2+  \leftrightarrow 7-$ & $2+  \leftrightarrow 7-$ & $8+  \leftrightarrow 8-$ & $8+  \leftrightarrow 8-$\\ 
$\hat 8$ & $2+  \leftrightarrow 7-$ & $1+  \leftrightarrow 1-$ & $2+  \leftrightarrow 7-$ & $1+  \leftrightarrow 1-$\\ 
 & & \\
$\hat 2$ & $3-  \leftrightarrow 3+$ & $3-  \leftrightarrow 3+$ & $1-  \leftrightarrow 4+$ & $1-  \leftrightarrow 4+$\\ 
$\hat 3$ & $1-  \leftrightarrow 4+$ & $2-  \leftrightarrow 2+$ & $1-  \leftrightarrow 4+$ & $2-  \leftrightarrow 2+$\\ 
 & & \\
$\hat 2$ & $5+  \leftrightarrow 3-$ & $5+  \leftrightarrow 3-$ & $4+  \leftrightarrow 4-$ & $4+  \leftrightarrow 4-$\\ 
$\hat 4$ & $1-  \leftrightarrow 3+$ & $2-  \leftrightarrow 2+$ & $1-  \leftrightarrow 3+$ & $2-  \leftrightarrow 2+$\\ 
 & & \\
$\hat 2$ & $5-  \leftrightarrow 5+$ & $5-  \leftrightarrow 5+$ & $4-  \leftrightarrow 6+$ & $4-  \leftrightarrow 6+$\\ 
$\hat 5$ & $1-  \leftrightarrow 3+$ & $2-  \leftrightarrow 2+$ & $1-  \leftrightarrow 3+$ & $2-  \leftrightarrow 2+$\\ 
 & & \\
$\hat 2$ & $7+  \leftrightarrow 5-$ & $7+  \leftrightarrow 5-$ & $6+  \leftrightarrow 6-$ & $6+  \leftrightarrow 6-$\\ 
$\hat 6$ & $1-  \leftrightarrow 3+$ & $2-  \leftrightarrow 2+$ & $1-  \leftrightarrow 3+$ & $2-  \leftrightarrow 2+$\\ 
 & & \\
$\hat 2$ & $7-  \leftrightarrow 7+$ & $7-  \leftrightarrow 7+$ & $6-  \leftrightarrow 8+$ & $6-  \leftrightarrow 8+$\\ 
$\hat 7$ & $1-  \leftrightarrow 3+$ & $2-  \leftrightarrow 2+$ & $1-  \leftrightarrow 3+$ & $2-  \leftrightarrow 2+$\\ 
 & & \\
$\hat 2$ & $1+  \leftrightarrow 7-$ & $1+  \leftrightarrow 7-$ & $8+  \leftrightarrow 8-$ & $8+  \leftrightarrow 8-$\\ 
$\hat 8$ & $1-  \leftrightarrow 3+$ & $2-  \leftrightarrow 2+$ & $1-  \leftrightarrow 3+$ & $2-  \leftrightarrow 2+$\\ 
 & & \\
$\hat 3$ & $4-  \leftrightarrow 4+$ & $4-  \leftrightarrow 4+$ & $2-  \leftrightarrow 5+$ & $2-  \leftrightarrow 5+$\\ 
$\hat 4$ & $2-  \leftrightarrow 5+$ & $3-  \leftrightarrow 3+$ & $2-  \leftrightarrow 5+$ & $3-  \leftrightarrow 3+$\\ 
\end{tabular}
\caption{\label{firstel} Elementary changes} 
\end{figure}

\begin{figure}
\begin{tabular}{c c c c c}
$\hat 3$ & $6+  \leftrightarrow 4-$ & $6+  \leftrightarrow 4-$ & $5+  \leftrightarrow 5-$ & $5+  \leftrightarrow 5-$\\ 
$\hat 5$ & $2-  \leftrightarrow 4+$ & $3-  \leftrightarrow 3+$ & $2-  \leftrightarrow 4+$ & $3-  \leftrightarrow 3+$\\ 
 & & \\
$\hat 3$ & $6-  \leftrightarrow 6+$ & $6-  \leftrightarrow 6+$ & $5-  \leftrightarrow 7+$ & $5-  \leftrightarrow 7+$\\ 
$\hat 6$ & $2-  \leftrightarrow 4+$ & $3-  \leftrightarrow 3+$ & $2-  \leftrightarrow 4+$ & $3-  \leftrightarrow 3+$\\ 
 & & \\
$\hat 3$ & $8+  \leftrightarrow 6-$ & $8+  \leftrightarrow 6-$ & $7+  \leftrightarrow 7-$ & $7+  \leftrightarrow 7-$\\ 
$\hat 7$ & $2-  \leftrightarrow 4+$ & $3-  \leftrightarrow 3+$ & $2-  \leftrightarrow 4+$ & $3-  \leftrightarrow 3+$\\ 
 & & \\
$\hat 3$ & $8-  \leftrightarrow 8+$ & $8-  \leftrightarrow 8+$ & $7-  \leftrightarrow 1+$ & $7-  \leftrightarrow 1+$\\ 
$\hat 8$ & $2-  \leftrightarrow 4+$ & $3-  \leftrightarrow 3+$ & $2-  \leftrightarrow 4+$ & $3-  \leftrightarrow 3+$\\ 
 & & \\
$\hat 4$ & $5-  \leftrightarrow 5+$ & $5-  \leftrightarrow 5+$ & $3-  \leftrightarrow 6+$ & $3-  \leftrightarrow 6+$\\ 
$\hat 5$ & $3-  \leftrightarrow 6+$ & $4-  \leftrightarrow 4+$ & $3-  \leftrightarrow 6+$ & $4-  \leftrightarrow 4+$\\ 
 & & \\
$\hat 4$ & $7+  \leftrightarrow 5-$ & $7+  \leftrightarrow 5-$ & $6+  \leftrightarrow 6-$ & $6+  \leftrightarrow 6-$\\ 
$\hat 6$ & $3-  \leftrightarrow 5+$ & $4-  \leftrightarrow 4+$ & $3-  \leftrightarrow 5+$ & $4-  \leftrightarrow 4+$\\ 
 & & \\
$\hat 4$ & $7-  \leftrightarrow 7+$ & $7-  \leftrightarrow 7+$ & $6-  \leftrightarrow 8+$ & $6-  \leftrightarrow 8+$\\ 
$\hat 7$ & $3-  \leftrightarrow 5+$ & $4-  \leftrightarrow 4+$ & $3-  \leftrightarrow 5+$ & $4-  \leftrightarrow 4+$\\ 
 & & \\
$\hat 4$ & $1+  \leftrightarrow 7-$ & $1+  \leftrightarrow 7-$ & $8+  \leftrightarrow 8-$ & $8+  \leftrightarrow 8-$\\ 
$\hat 8$ & $3-  \leftrightarrow 5+$ & $4-  \leftrightarrow 4+$ & $3-  \leftrightarrow 5+$ & $4-  \leftrightarrow 4+$\\ 
 & & \\
$\hat 5$ & $6-  \leftrightarrow 6+$ & $6-  \leftrightarrow 6+$ & $4-  \leftrightarrow 7+$ & $4-  \leftrightarrow 7+$\\ 
$\hat 6$ & $4-  \leftrightarrow 7+$ & $5-  \leftrightarrow 5+$ & $4-  \leftrightarrow 7+$ & $5-  \leftrightarrow 5+$\\ 
 & & \\
$\hat 5$ & $8+  \leftrightarrow 6-$ & $8+  \leftrightarrow 6-$ & $7+  \leftrightarrow 7-$ & $7+  \leftrightarrow 7-$\\ 
$\hat 7$ & $4-  \leftrightarrow 6+$ & $5-  \leftrightarrow 5+$ & $4-  \leftrightarrow 6+$ & $5-  \leftrightarrow 5+$\\ 
 & & \\
$\hat 5$ & $8-  \leftrightarrow 8+$ & $8-  \leftrightarrow 8+$ & $7-  \leftrightarrow 1+$ & $7-  \leftrightarrow 1+$\\ 
$\hat 8$ & $4-  \leftrightarrow 6+$ & $5-  \leftrightarrow 5+$ & $4-  \leftrightarrow 6+$ & $5-  \leftrightarrow 5+$\\ 
 & & \\
$\hat 6$ & $7-  \leftrightarrow 7+$ & $7-  \leftrightarrow 7+$ & $5-  \leftrightarrow 8+$ & $5-  \leftrightarrow 8+$\\ 
$\hat 7$ & $5-  \leftrightarrow 8+$ & $6-  \leftrightarrow 6+$ & $5-  \leftrightarrow 8+$ & $6-  \leftrightarrow 6+$\\ 
 & & \\
$\hat 6$ & $1+  \leftrightarrow 7-$ & $1+  \leftrightarrow 7-$ & $8+  \leftrightarrow 8-$ & $8+  \leftrightarrow 8-$\\ 
$\hat 8$ & $5-  \leftrightarrow 7+$ & $6-  \leftrightarrow 6+$ & $5-  \leftrightarrow 7+$ & $6-  \leftrightarrow 6+$\\ 
 & & \\
$\hat 7$ & $8-  \leftrightarrow 8+$ & $8-  \leftrightarrow 8+$ & $6-  \leftrightarrow 1+$ & $6-  \leftrightarrow 1+$\\ 
$\hat 8$ & $6-  \leftrightarrow 1+$ & $7-  \leftrightarrow 7+$ & $6-  \leftrightarrow 1+$ & $7-  \leftrightarrow 7+$\\ 
\end{tabular}
\caption{\label{secondel} Elementary changes, continued} 
\end{figure}

\vspace{2cm}

severine.fiedler@live.fr

\end{document}